\documentclass[11pt]{article}
\usepackage[]{amsmath,amssymb}

\newtheorem{theorem}{Theorem}[section]
\newtheorem{lemma}[theorem]{Lemma}
\newtheorem{proposition}[theorem]{Proposition}

\newtheorem{exAux}[theorem]{Example}

\newtheorem{Def}[theorem]{Definition}
\newenvironment{definition}{\begin{Def} \rm}{\end{Def}}
\newtheorem{Note}[theorem]{Note}
\newenvironment{note}{\begin{Note} \rm}{\end{Note}}
\newtheorem{Problem}[theorem]{Problem}

\newtheorem{Rem}[theorem]{Remark}
\newenvironment{remark}{\begin{Rem} \rm}{\end{Rem}}
\newtheorem{Not}[theorem]{Notation}
\newenvironment{notation}{\begin{Not} \rm}{\end{Not}}
\newtheorem{Ass}[theorem]{Assumption}

\newenvironment{proof}{\medskip\noindent{\bf Proof.\ }}{\qed\medskip}
\newenvironment{proofof}[1]{\medskip\noindent{\bf Proof  of {#1}.\ 
}}{\qed\medskip}
\newcommand{\qed}{\hfill\mbox{$\Box$\qquad\qquad}}

\newcommand{\Mat}[1]{\text{\rm Mat}_{#1}(\mathbb{K})}

\renewcommand{\indent}{\hspace{6mm}}

\begin{document}
\thispagestyle{empty}

\begin{center}
\LARGE \bf
\noindent
Affine transformations of a Leonard pair
\end{center}

\smallskip

\begin{center}
\Large
Kazumasa Nomura and Paul Terwilliger
\end{center}

\smallskip

\begin{quote}
\small 
\begin{center}
\bf Abstract
\end{center}

\indent
Let $\mathbb{K}$ denote a field and let $V$ denote a vector space
over $\mathbb{K}$ with finite positive dimension. 
We consider an ordered pair of linear transformations 
$A : V \to V$ and $A^* : V \to V$ that satisfy (i) and (ii) below:
\begin{itemize}
\item[(i)] There exists a basis for $V$ with respect to which the
matrix representing $A$ is irreducible tridiagonal and the matrix
representing $A^*$ is diagonal.
\item[(ii)] There exists a basis for $V$ with respect to which the
matrix representing $A^*$ is irreducible tridiagonal and the matrix
representing $A$ is diagonal.
\end{itemize}
We call such a pair a {\em Leonard pair} on $V$.
Let $\xi, \zeta, \xi^*, \zeta^*$ denote scalars in $\mathbb{K}$ with
$\xi, \xi^*$ nonzero,
and note that $\xi A+\zeta I$, $\xi^* A^* + \zeta^* I$ is a Leonard 
pair on $V$.
We give necessary and sufficient conditions for this Leonard pair
to be isomorphic to $A,A^*$. We also give necessary and sufficient
conditions for this Leonard pair to be isomorphic to the
Leonard pair $A^*,A$.
\end{quote}

\section{Leonard pairs}

\indent
We begin by recalling the notion of a Leonard pair.
We will use the following terms.
A square matrix $X$ is said to be {\em tridiagonal}
whenever each nonzero entry lies on either the diagonal, the subdiagonal,
or the superdiagonal. Assume $X$ is tridiagonal.
Then $X$ is said to be {\em irreducible}
whenever each entry on the subdiagonal is nonzero and each entry on
the superdiagonal is nonzero.
We now define a Leonard pair.
For the rest of this paper $\mathbb{K}$ will denote a field.

\begin{definition}  \cite{T:Leonard}  \label{def:LP}      \samepage
Let $V$ denote a vector space over $\mathbb{K}$ with finite positive
dimension.
By a {\em Leonard pair} on $V$ we mean an ordered pair $A,A^*$
where $A:V \to V$ and $A^*:V \to V$ are linear transformations
that satisfy (i) and (ii) below:
\begin{itemize}
\item[(i)] There exists a basis for $V$ with respect to which the
matrix representing $A$ is irreducible tridiagonal and the matrix
representing $A^*$ is diagonal.
\item[(ii)] There exists a basis for $V$ with respect to which the
matrix representing $A^*$ is irreducible tridiagonal and the matrix
representing $A$ is diagonal.
\end{itemize}
\end{definition}

\begin{note}
It is a common notational convention to use $A^*$ to represent the
conjugate-transpose of $A$. We are not using this convention.
In a Leonard pair $A,A^*$ the linear transformations $A$ and
$A^*$ are arbitrary subject to (i) and (ii) above.
\end{note}

We refer the reader to 
\cite{H},
\cite{N:aw},
\cite{NT:balanced}, \cite{NT:formula}, \cite{NT:det}, \cite{NT:mu},
\cite{NT:span}, \cite{NT:switch},
\cite{P}, \cite{T:sub1}, \cite{T:sub3}, \cite{T:Leonard},
\cite{T:24points}, \cite{T:canform}, \cite{T:intro},
\cite{T:intro2}, \cite{T:split}, \cite{T:array}, \cite{T:qRacah},
\cite{T:survey}, \cite{TV}, \cite{V}, \cite{V:AW}
for background on Leonard pairs.
We especially recommend the survey \cite{T:survey}.
See \cite{AC}, \cite{AC2}, \cite{BT:Borel}, \cite{BT:loop},
\cite{H:tetra}, \cite{HT:tetra},
\cite{ITT}, \cite{IT:shape},
\cite{IT:uqsl2hat}, \cite{IT:non-nilpotent}, 
\cite{IT:tetra}, \cite{IT:inverting},
\cite{ITW:equitable}, 
\cite{N:refine}, \cite{N:height1},
\cite{R:multi}, \cite{R:6j},
\cite{T:qSerre}, \cite{T:Kac-Moody},
\cite{Z}
for related topics.

In this paper we consider the following situation.
Let $V$ denote a vector space over $\mathbb{K}$ with finite
positive dimension and let $A,A^*$ denote
a Leonard pair on $V$.  Let
$\xi,\zeta,\xi^*,\zeta^*$ denote scalars in $\mathbb{K}$ with 
$\xi,\xi^*$ nonzero,
and note that $\xi A+\zeta I$, $\xi^* A^* + \zeta^* I$ is a Leonard pair
on $V$. We give necessary and sufficient conditions
for this Leonard pair to be isomorphic to $A,A^*$.
We also give necessary and sufficient conditions
for this Leonard pair to be isomorphic to the
Leonard pair $A^*,A$.

\section{Leonard systems}

\indent
When working with a Leonard pair, it is convenient to consider a closely
related object called a {\em Leonard system}. 
To prepare for our definition
of a Leonard system, we recall a few concepts from linear algebra.
Let $d$ denote a nonnegative integer and let
$\Mat{d+1}$ denote the $\mathbb{K}$-algebra consisting of all $d+1$ by
$d+1$ matrices that have entries in $\mathbb{K}$. 
We index the rows and  columns by $0, 1, \ldots, d$. 
We let $\mathbb{K}^{d+1}$ denote the $\mathbb{K}$-vector space of all
$d+1$ by $1$ matrices that have entries in $\mathbb{K}$. We index the
rows by $0,1,\ldots,d$. We view $\mathbb{K}^{d+1}$ as a left module
for $\Mat{d+1}$. We observe this module is irreducible. 
For the rest of this paper, let $\cal A$ denote a $\mathbb{K}$-algebra 
isomorphic to $\Mat{d+1}$ and let $V$ denote an irreducible left $\cal A$-module.
We remark that $V$ is unique
up to isomorphism of $\cal A$-modules, and that $V$ has dimension $d+1$.
Let $\{v_i\}_{i=0}^d$ denote a basis for $V$.
For $X \in {\cal A}$ and $Y \in \Mat{d+1}$, we say 
{\em $Y$ represents $X$ with respect to} $\{v_i\}_{i=0}^d$
whenever $X v_j = \sum_{i=0}^d Y_{ij}v_i$ for $0 \leq j \leq d$.
For $A \in \cal A$ we say $A$ is {\em multiplicity-free}
whenever it has $d+1$ mutually distinct eigenvalues in $\mathbb{K}$. 
Assume $A$ is multiplicity-free. 
Let $\{\theta_i\}_{i=0}^d$ denote an ordering 
of the eigenvalues of $A$, and for $0 \leq i \leq d$ put
\begin{equation}        \label{eq:defEi}
    E_i = \prod_{\stackrel{0 \leq j \leq d}{j\neq i}}
             \frac{A-\theta_j I}{\theta_i - \theta_j},
\end{equation}
where $I$ denotes the identity of $\cal A$. 
We observe
(i) $AE_i = \theta_i E_i$ $(0 \leq i \leq d)$;
(ii) $E_i E_j = \delta_{i,j} E_i$ $(0 \leq i,j \leq d)$;
(iii) $\sum_{i=0}^{d} E_i = I$;
(iv) $A = \sum_{i=0}^{d} \theta_i E_i$.
Let $\cal D$ denote the subalgebra of $\cal A$ generated by $A$.
Using (i)--(iv) we find the sequence $\{E_i\}_{i=0}^d$
is a basis for the $\mathbb{K}$-vector space $\cal D$.
We call $E_i$ the {\em primitive idempotent} of $A$ associated with
$\theta_i$. 
It is helpful to think of these primitive idempotents as follows.
Observe 
\[
 V=E_0V+E_1V+\cdots+E_dV  \qquad \text{(direct sum)}.
\]
For $0 \leq i \leq d$, $E_iV$ is the (one dimensional) eigenspace of $A$
in $V$ associated with the eigenvalue $\theta_i$, and $E_i$ acts on $V$
as the projection onto this eigenspace.

\medskip

By a {\em Leonard pair in $\cal A$} we mean an ordered pair of elements
taken from $\cal A$ that act on $V$ as a Leonard pair in the sense of
Definition \ref{def:LP}.
We call $\cal A$ the {\em ambient algebra} of the pair and say the
pair is {\em over $\mathbb{K}$}.
We now define a Leonard system.

\begin{definition}  \cite{T:Leonard}     \label{def:LS}   \samepage
By a {\em Leonard system} in $\cal A$ we mean a sequence
\[
  \Phi= (A; \{E_i\}_{i=0}^d; A^*; \{E^*_i\}_{i=0}^d)
\]
that satisfies (i)--(v) below.
\begin{itemize}
\item[(i)] Each of $A$, $A^*$ is a multiplicity-free element in $\cal A$.
\item[(ii)] $\{E_i\}_{i=0}^d$  is an ordering of the
   primitive idempotents of $A$.
\item[(iii)] $\{E^*_i\}_{i=0}^d$ is an ordering of the
   primitive idempotents of $A^*$.
\item[(iv)] For $0 \leq i,j \leq d$, 
\begin{equation}           \label{eq:Astrid}
   E_i A^* E_j =
    \begin{cases}  
        0 & \text{\rm if $|i-j|>1$},  \\
        \neq 0 & \text{\rm if $|i-j|=1$}.
    \end{cases}
\end{equation}
\item[(v)] For $0 \leq i,j \leq d$, 
\begin{equation}             \label{eq:Atrid}
   E^*_i A E^*_j =
    \begin{cases}  
        0 & \text{\rm if $|i-j|>1$},  \\
        \neq 0 & \text{\rm if $|i-j|=1$}.
    \end{cases}
\end{equation}
\end{itemize}
We refer to $d$ as the {\em diameter} of $\Phi$ and
say {\em $\Phi$ is over $\mathbb{K}$}.
We call $\cal A$ the {\em ambient algebra} of $\Phi$.
\end{definition}

\medskip

Leonard systems are related to Leonard pairs as follows.
Let $(A; \{E_i\}_{i=0}^d;A^*; \{E^*_i\}_{i=0}^d)$ denote a Leonard system
in $\cal A$. Then $A,A^*$ is a Leonard pair in $\cal A$
\cite[Section 3]{T:qRacah}.
Conversely, suppose $A,A^*$ is a Leonard pair in $\cal A$.
Then each of $A,A^*$ is multiplicity-free \cite[Lemma 1.3]{T:Leonard}.
Moreover there exists an ordering $\{E_i\}_{i=0}^d$ of the
primitive idempotents of $A$, and 
there exists an ordering $\{E^*_i\}_{i=0}^d$  of the
primitive idempotents of $A^*$, such that
$(A; \{E_i\}_{i=0}^d; A^*; \{E^*_i\}_{i=0}^d)$
is a Leonard system in $\cal A$ \cite[Lemma 3.3]{T:qRacah}.
We say this Leonard system is {\em associated} with
the Leonard pair $A,A^*$.

\medskip

We recall the notion of {\em isomorphism} for Leonard pairs
and Leonard systems.

\begin{definition}            \label{def:LPiso}    \samepage
Let $A,A^*$ and $B,B^*$ denote Leonard pairs over $\mathbb{K}$.
By an {\em isomorphism of Leonard pairs from $A,A^*$ to $B,B^*$}
we mean an isomorphism of $\mathbb{K}$-algebras
from the ambient algebra of $A,A^*$ to the ambient algebra $B,B^*$
that sends $A$ to $B$ and $A^*$ to $B^*$.
The Leonard pairs $A,A^*$ and $B,B^*$ are said to be {\em isomorphic}
whenever there exists an isomorphism of Leonard pairs from
$A,A^*$ to $B,B^*$.
\end{definition}

\medskip

Let $\Phi$ denote the Leonard system from Definition \ref{def:LS} and
let $\sigma:{\cal A} \to {\cal A}'$ denote an isomorphism of
$\mathbb{K}$-algebras.
We write
$\Phi^\sigma 
  :=(A^\sigma; \{E^\sigma_i\}_{i=0}^d;
     A^{*\sigma}; \{E^{*\sigma}_i\}_{i=0}^d)$
and observe $\Phi^\sigma$ is a Leonard system in ${\cal A}'$.

\medskip

\begin{definition}              \label{def:LSiso}  \samepage
Let $\Phi$ and $\Phi'$ denote Leonard systems over $\mathbb{K}$.
By an {\em isomorphism of Leonard systems from $\Phi$ to $\Phi'$}
we mean an isomorphism of $\mathbb{K}$-algebras $\sigma$
from the ambient algebra of $\Phi$ to the ambient algebra of $\Phi'$
such that $\Phi^\sigma=\Phi'$.
The Leonard systems $\Phi$ and $\Phi'$ are said to be {\em isomorphic} whenever
there exists an isomorphism of Leonard systems from $\Phi$ to $\Phi'$.
\end{definition}

\section{The $D_4$ action}

\indent
Let $\Phi=(A; \{E_i\}_{i=0}^d; A^*; \{E^*_i\}_{i=0}^d)$
denote a Leonard system in $\cal A$.
Then each of the following is a Leonard system in $\cal A$:
\begin{eqnarray*}
\Phi^{*}  &:=& 
       (A^*; \{E^*_i\}_{i=0}^d; A; \{E_i\}_{i=0}^d), \\
\Phi^{\downarrow} &:=&
       (A; \{E_i\}_{i=0}^d; A^*; \{E^*_{d-i}\}_{i=0}^d), \\
\Phi^{\Downarrow} &:=&
       (A; \{E_{d-i}\}_{i=0}^d; A^*; \{E^*_{i}\}_{i=0}^d).
\end{eqnarray*}
Viewing $*$, $\downarrow$, $\Downarrow$ as permutations on the set of
all the Leonard systems,
\begin{equation}    \label{eq:relation1}
\text{$*^2$$\;=\;$$\downarrow^2$$\;=\;$$\Downarrow^2$$\;=\;$$1$},
\end{equation}
\begin{equation}    \label{eq:relation2}
\text{$\Downarrow$$*$$\;=\;$$*$$\downarrow$}, \quad
\text{$\downarrow$$*$$\;=\;$$*$$\Downarrow$}, \quad
\text{$\downarrow$$\Downarrow$$\;=\;$$\Downarrow$$\downarrow$}.
\end{equation}
The group generated by symbols $*$, $\downarrow$, $\Downarrow$ subject
to the relations (\ref{eq:relation1}), (\ref{eq:relation2}) is the
dihedral group $D_4$. We recall $D_4$ is the group of symmetries of a
square, and has $8$ elements.
Apparently $*$, $\downarrow$, $\Downarrow$ induce an action of $D_4$
on the set of all Leonard systems.
Two Leonard systems will be called {\em relatives} whenever they are
in the same orbit of this $D_4$ action. 
The relatives of $\Phi$ are as follows:

\medskip

\begin{center}
\begin{tabular}{c|c}
name  &  relative \\
\hline
$\Phi$ & 
       $(A; \{E_i\}_{i=0}^d; A^*;  \{E^*_i\}_{i=0}^d)$ \\ 
$\Phi^{\downarrow}$ &
       $(A; \{E_i\}_{i=0}^d; A^*;  \{E^*_{d-i}\}_{i=0}^d)$ \\ 
$\Phi^{\Downarrow}$ &
       $(A; \{E_{d-i}\}_{i=0}^d; A^*;  \{E^*_i\}_{i=0}^d)$ \\ 
$\Phi^{\downarrow \Downarrow}$ &
       $(A; \{E_{d-i}\}_{i=0}^d; A^*;  \{E^*_{d-i}\}_{i=0}^d)$ \\ 
$\Phi^{*}$  & 
       $(A^*; \{E^*_i\}_{i=0}^d; A;  \{E_i\}_{i=0}^d)$ \\ 
$\Phi^{\downarrow *}$ &
       $(A^*; \{E^*_{d-i}\}_{i=0}^d; A;  \{E_i\}_{i=0}^d)$ \\ 
$\Phi^{\Downarrow *}$ &
       $(A^*; \{E^*_i\}_{i=0}^d; A;  \{E_{d-i}\}_{i=0}^d)$ \\ 
$\Phi^{\downarrow \Downarrow *}$ &
       $(A^*; \{E^*_{d-i}\}_{i=0}^d; A;  \{E_{d-i}\}_{i=0}^d)$
\end{tabular}
\end{center}

\section{The parameter array}

\indent
In this section we recall the parameter array of a Leonard system.

\begin{definition}        \label{def:th}       \samepage
Let $\Phi=(A; \{E_i\}_{i=0}^d;A^*; \{E^*_i\}_{i=0}^d)$ denote a 
Leonard system over $\mathbb{K}$.
For $0 \leq i \leq d$ we let $\theta_i$ (resp. $\theta^*_i$)
denote the eigenvalue of $A$ (resp. $A^*$) associated with
$E_i$ (resp. $E^*_i$).
We refer to $\{\theta_i\}_{i=0}^d$ (resp. $\{\theta^*_i\}_{i=0}^d$)
as the {\em eigenvalue sequence} (resp. {\em dual eigenvalue sequence})
of $\Phi$.
We observe $\{\theta_i\}_{i=0}^d$ (resp. $\{\theta^*_i\}_{i=0}^d$)
are mutually distinct
and contained in $\mathbb{K}$.
\end{definition}

\begin{definition}  \cite[Theorem 4.6]{NT:formula}  \label{def:splitseq} \samepage
Let $\Phi=(A; \{E_i\}_{i=0}^d;A^*; \{E^*_i\}_{i=0}^d)$ denote a 
Leonard system with eigenvalue sequence $\{\theta_i\}_{i=0}^d$
and dual eigenvalue sequence $\{\theta^*_i\}_{i=0}^d$.
For $1 \leq i \leq d$ we define
\begin{eqnarray}
 \varphi_i 
   &:=& (\theta^*_0-\theta^*_i)
       \frac{\text{tr}(E^*_0 \prod_{h=0}^{i-1}(A-\theta_hI))}
            {\text{tr}(E^*_0\prod_{h=0}^{i-2}(A-\theta_hI))},
                               \label{eq:defvarphi}   \\
 \phi_i
   &:=& (\theta^*_0-\theta^*_i)
       \frac{\text{tr}(E^*_0\prod_{h=0}^{i-1}(A-\theta_{d-h}I))}
            {\text{tr}(E^*_0\prod_{h=0}^{i-2}(A-\theta_{d-h}I))}
                               \label{eq:defphi},
\end{eqnarray}
where tr means trace.
In (\ref{eq:defvarphi}), (\ref{eq:defphi}) the denominators 
are nonzero by \cite[Corollary 4.5]{NT:formula}.
The sequence $\{\varphi_i\}_{i=1}^d$ 
(resp. $\{\phi_i\}_{i=1}^d$) is called the 
{\em first split sequence} (resp. {\em second split sequence}) of $\Phi$.
\end{definition}

\begin{definition}            \label{def:param}        \samepage
Let $\Phi=(A; \{E_i\}_{i=0}^d; A^*; \{E^*_i\}_{i=0}^d)$
denote a Leonard system over $\mathbb{K}$.
By the {\em parameter array of $\Phi$} we mean the sequence
$(\{\theta_i\}_{i=0}^d; \{\theta^*_i\}_{i=0}^d;
        \{\varphi_i\}_{i=1}^d; \{\phi_i\}_{i=1}^d)$,
where the $\theta_i$, $\theta^*_i$ are from
Definition \ref{def:th} and the $\varphi_i$, $\phi_i$ are
from Definition \ref{def:splitseq}.
\end{definition}

\begin{theorem}  \cite[Theorem 1.9]{T:Leonard} \label{thm:classify} \samepage
Let $d$ denote a nonnegative integer and let
\begin{equation}            \label{eq:paramarray}
(\{\theta_i\}_{i=0}^d; \{\theta^*_i\}_{i=0}^d;
        \{\varphi_i\}_{i=1}^d; \{\phi_i\}_{i=1}^d)
\end{equation}
denote a sequence of scalars taken from $\mathbb{K}$.
Then there exists a Leonard system $\Phi$ over $\mathbb{K}$
with parameter array (\ref{eq:paramarray}) if and only if
(PA1)--(PA5) hold below.
\begin{itemize}
\item[]
\begin{itemize}
\item[(PA1)]  $\varphi_i \neq 0$, $\phi_i \neq 0$ $(1 \leq i \leq d)$.
\item[(PA2)]  $\theta_i \neq \theta_j$, $\theta^*_i \neq \theta^*_j$
   if $i \neq j$ $(0 \leq i,j \leq d)$.
\item[(PA3)] For $1 \leq i \leq d$,
\[
 \varphi_i = \phi_1 
 \sum_{h=0}^{i-1} \frac{\theta_h - \theta_{d-h}}
                       {\theta_0 - \theta_d}
         + (\theta^*_{i}-\theta^*_{0})(\theta_{i-1}-\theta_{d}).
\]
\item[(PA4)] For $1 \leq i \leq d$,
\[
\phi_i = \varphi_1
\sum_{h=0}^{i-1} \frac{\theta_h - \theta_{d-h}}
                       {\theta_0 - \theta_d}
         + (\theta^*_{i}-\theta^*_{0})(\theta_{d-i+1}-\theta_{0}).
\]
\item[(PA5)]  The expressions
\begin{equation}        \label{eq:indep}
   \frac{\theta_{i-2}-\theta_{i+1}}{\theta_{i-1}-\theta_{i}},
 \quad
   \frac{\theta^*_{i-2}-\theta^*_{i+1}}{\theta^*_{i-1}-\theta^*_{i}}
\end{equation}
are equal and independent of $\;i\;$ for $2 \leq i \leq d-1$.
\end{itemize}
\end{itemize}
Suppose (PA1)--(PA5) hold. Then $\Phi$ is unique up to isomorphism
of Leonard systems.
\end{theorem}

\medskip

The $D_4$ action affects the parameter array as follows.

\medskip

\begin{lemma}   \cite[Theorem 1.11]{T:Leonard}        \label{lem:D4}   \samepage
Let $\Phi=(A; \{E_i\}_{i=0}^d; A^*; \{E^*_i\}_{i=0}^d)$
denote a Leonard system with parameter array
 $(\{\theta_i\}_{i=0}^d; \{\theta^*_{i}\}_{i=0}^d;
        \{\varphi_{i}\}_{i=1}^d; \{\phi_{i}\}_{i=1}^d)$.
For each relative of $\Phi$ the parameter array is given below.
\[
  \begin{array}{c|c}
\rule[-1ex]{0cm}{3ex}   \text{\rm relative} & \text{\rm parameter array} \\
\hline
\rule{0cm}{3ex} \Phi &
  (\{\theta_i\}_{i=0}^d; \{\theta^*_{i}\}_{i=0}^d;
        \{\varphi_{i}\}_{i=1}^d; \{\phi_{i}\}_{i=1}^d)
\\
\rule{0cm}{3ex} \Phi^{\downarrow} &
   (\{\theta_i\}_{i=0}^d; \{\theta^*_{d-i}\}_{i=0}^d;
        \{\phi_{d-i+1}\}_{i=1}^d; \{\varphi_{d-i+1}\}_{i=1}^d)  
\\
\rule{0cm}{3ex} \Phi^{\Downarrow} &
  (\{\theta_{d-i}\}_{i=0}^d; \{\theta^*_{i}\}_{i=0}^d;
        \{\phi_{i}\}_{i=1}^d; \{\varphi_{i}\}_{i=1}^d)
\\
\rule{0cm}{3ex} \Phi^{\downarrow\Downarrow} &
  (\{\theta_{d-i}\}_{i=0}^d; \{\theta^*_{d-i}\}_{i=0}^d;
        \{\varphi_{d-i+1}\}_{i=1}^d; \{\phi_{d-i+1}\}_{i=1}^d)
\\
\rule{0cm}{3ex} \Phi^*  &
   (\{\theta^*_i\}_{i=0}^d; \{\theta_i\}_{i=0}^d;
        \{\varphi_i\}_{i=1}^d; \{\phi_{d-i+1}\}_{i=1}^d)
\\
\rule{0cm}{3ex} \rule{0cm}{3ex} \Phi^{\downarrow*} &
   (\{\theta^*_{d-i}\}_{i=0}^d; \{\theta_i\}_{i=0}^d;
        \{\phi_{d-i+1}\}_{i=1}^d; \{\varphi_{i}\}_{i=1}^d)
\\
\rule{0cm}{3ex} \Phi^{\Downarrow*} &
  (\{\theta^*_{i}\}_{i=0}^d; \{\theta_{d-i}\}_{i=0}^d;
        \{\phi_{i}\}_{i=1}^d; \{\varphi_{d-i+1}\}_{i=1}^d)
\\
\rule{0cm}{3ex} \Phi^{\downarrow\Downarrow*} &
   (\{\theta^*_{d-i}\}_{i=0}^d; \{\theta_{d-i}\}_{i=0}^d;
        \{\varphi_{d-i+1}\}_{i=1}^d; \{\phi_{i}\}_{i=1}^d)
\end{array}
\]
\end{lemma}

\section{Affine transformations of a Leonard system}

In this section we consider the affine transformations of
a Leonard system. We start with an observation.

\begin{lemma}          \label{lem:affine}    \samepage
Let $\Phi=(A; \{E_i\}_{i=0}^d;A^*; \{E^*_i\}_{i=0}^d)$ denote a 
Leonard system in $\cal A$.
Let $\xi, \zeta, \xi^*, \zeta^*$ denote scalars in $\mathbb{K}$
with $\xi, \xi^*$ nonzero.
Then the sequence
\begin{equation}         \label{eq:affPhi}
     (\xi A + \zeta I; \{E_i\}_{i=0}^d; 
             \xi^* A^* +\zeta^* I; \{E^*_i\}_{i=0}^d)
\end{equation}
is a Leonard system in $\cal A$.
\end{lemma}

\begin{definition}            \label{def:aff}        \samepage
Referring to Lemma \ref{lem:affine},
we call (\ref{eq:affPhi}) the {\em affine transformation}
of $\Phi$ associated with $\xi, \zeta, \xi^*, \zeta^*$.
\end{definition}

\begin{definition}           \label{def:affiso}        \samepage
Let $\Phi$ and $\Phi'$ denote Leonard systems over $\mathbb{K}$.
We say $\Phi$ and $\Phi'$ are {\em affine isomorphic} whenever
$\Phi$ is isomorphic to an affine transformation of $\Phi'$.
Observe that affine isomorphism is an equivalence relation.
\end{definition}

\medskip

Let $\Phi$ denote a Leonard system. We now consider how the
set of relatives of $\Phi$ is partitioned into affine isomorphism
classes.
In order to avoid trivialities we assume the diameter of $\Phi$
is at least $1$.
The following is our main result on this topic.

\begin{theorem}        \label{thm:main}        \samepage
Let $\Phi$ denote a Leonard system with first split sequence
$\{\varphi_i\}_{i=1}^d$ and second split sequence $\{\phi_i\}_{i=1}^d$.
Assume $d \geq 1$.
\begin{itemize}
\item[(i)] 
Assume $\varphi_1=\varphi_d=-\phi_1=-\phi_d$.
Then all eight relatives of $\Phi$ are mutually affine isomorphic.
\item[(ii)]
Assume $\varphi_1=\varphi_d$, $\phi_1=\phi_d$ and $\varphi_1 \neq -\phi_1$.
Then the relatives of $\Phi$ form exactly two affine isomorphism classes,
consisting of
 $\{\Phi, \Phi^{\downarrow\Downarrow}, \Phi^*,  \Phi^{\downarrow\Downarrow*}\}$,
 $\{\Phi^{\downarrow},\Phi^{\Downarrow},\Phi^{\downarrow*},\Phi^{\Downarrow*}\}$.
\item[(iii)]
Assume $\varphi_1=\varphi_d$ and $\phi_1 \neq \phi_d$.
Then the relatives of $\Phi$ form  exactly four affine isomorphism classes,
consisting of
 $\{\Phi, \Phi^{\downarrow\Downarrow*}\}$,
 $\{\Phi^{\downarrow}, \Phi^{\downarrow*}\}$, 
 $\{\Phi^{\Downarrow}, \Phi^{\Downarrow*} \}$,
 $\{\Phi^{\downarrow\Downarrow}, \Phi^{*} \}$.
\item[(iv)] 
Assume $\phi_1 = \phi_d$ and $\varphi_1 \neq \varphi_d$.
Then the relatives of $\Phi$ form exactly four affine isomorphism classes,
consisting of
 $\{\Phi, \Phi^*\}$,
 $\{\Phi^{\downarrow}, \Phi^{\Downarrow*}\}$,
 $\{\Phi^{\Downarrow}, \Phi^{\downarrow*} \}$,
 $\{\Phi^{\downarrow\Downarrow}, \Phi^{\downarrow\Downarrow*} \}$.
\item[(v)]
Assume $\varphi_1=-\phi_1$, $\varphi_d=-\phi_d$ and $\varphi_1 \neq \varphi_d$.
Then the relatives of $\Phi$ form exactly four affine isomorphism classes,
consisting of
 $\{\Phi, \Phi^{\Downarrow}\}$,
 $\{\Phi^{\downarrow}, \Phi^{\downarrow\Downarrow} \}$,
 $\{\Phi^{*}, \Phi^{\Downarrow*} \}$,
 $\{\Phi^{\downarrow*}, \Phi^{\downarrow\Downarrow*} \}$.
\item[(vi)]
Assume $\varphi_1=-\phi_d$, $\varphi_d=-\phi_1$ and $\varphi_1 \neq \varphi_d$.
Then the relatives of $\Phi$ form exactly four affine isomorphism classes,
consisting of
 $\{\Phi, \Phi^{\downarrow}\}$,
 $\{\Phi^{\Downarrow}, \Phi^{\downarrow\Downarrow}\}$,
 $\{\Phi^{*}, \Phi^{\downarrow*}\}$,
 $\{\Phi^{\Downarrow*}, \Phi^{\downarrow\Downarrow*} \}$.
\item[(vii)]
Assume none of (i)--(vi) hold above.
Then $\varphi_1 \neq \varphi_d$,
$\phi_1 \neq \phi_d$,
at least one of $\varphi_1 \neq  -\phi_1$,
$\varphi_d \neq -\phi_d$,
and at least one of $\varphi_1 \neq -\phi_d$,
$\varphi_d \neq -\phi_1$.
In this case the eight relatives of $\Phi$ are mutually non affine
isomorphic.
\end{itemize}
\end{theorem}

\medskip

The proof of Theorem \ref{thm:main} will be given in Section 9. 
In Sections 6--8 we obtain some results that will
be used in this proof.

\section{How the parameter array is affected by affine transformation}

\indent
Let $\Phi$ denote a Leonard system.
In this section we consider 
how the parameter array of $\Phi$ is affected by affine transformation.

\begin{lemma}          \label{lem:affparam}     \samepage
Referring to Lemma \ref{lem:affine},
let 
 $(\{\theta_i\}_{i=0}^d; \{\theta^*_i\}_{i=0}^d;
    \{\varphi_i\}_{i=1}^d; \{\phi_i\}_{i=1}^d)$
denote the parameter array of $\Phi$.
Then the parameter array of the Leonard system (\ref{eq:affPhi}) is
\begin{equation}        \label{eq:affparam}
 (\{\xi \theta_i+\zeta\}_{i=0}^d;
  \{\xi^* \theta^*_i + \zeta^*\}_{i=0}^d;
  \{\xi\xi^* \varphi_i\}_{i=1}^d;
  \{\xi\xi^* \phi_i\}_{i=1}^d).
\end{equation}
\end{lemma}

\begin{proof}
By Definition \ref{def:th}, for 0 $\leq i \leq d$ the scalar
$\theta_i$ is the eigenvalue of $A$ associated with $E_i$, 
so $\xi \theta_i+ \zeta$ is the eigenvalue of $\xi A + \zeta I$
associated with $E_i$. Thus $\{\xi \theta_i+\zeta\}_{i=0}^d$ 
is the eigenvalue sequence of (\ref{eq:affPhi}). 
Similarly $\{\xi^* \theta^*_i+\zeta^*\}_{i=0}^d$ is the dual eigenvalue 
sequence of (\ref{eq:affPhi}). 
In the right-hand side of (\ref{eq:defvarphi}), if we replace $A$ by
$\xi A+\zeta I$, and if we replace $\theta_j$, $\theta^*_j$
by $\xi\theta_j+\zeta$, $\xi^* \theta^*_j + \zeta^*$ $(0 \leq j \leq d)$ 
and simplify the result we get $\xi \xi^* \varphi_i$. 
Therefore $\{\xi \xi^* \varphi_i\}_{i=1}^d$ is
the first split sequence of (\ref{eq:affPhi}). 
Similarly
$\{\xi \xi^* \phi_i\}_{i=1}^d$ is the second split sequence of 
(\ref{eq:affPhi}) and the result follows.
\end{proof}

\section{Some equations}

\indent
In this section we obtain some equations that will
be useful in the proof of Theorem \ref{thm:main}.

\begin{notation}         \label{notation}
Let $\Phi=(A; \{E_i\}_{i=0}^d; A^*; \{E^*_i\}_{i=0}^d)$ 
denote a Leonard system over $\mathbb{K}$, 
with parameter array 
$(\{\theta_i\}_{i=0}^d; \{\theta^*_i\}_{i=0}^d$;
  $\{\varphi_i\}_{i=1}^d; \{\phi_i\}_{i=1}^d)$.
To avoid trivialities we assume $d\geq 1$.
\end{notation}

\begin{lemma}    \cite[Lemma 9.5]{T:Leonard}  \label{lem:thh} \samepage
Referring to Notation \ref{notation},
\[
  \frac{\theta_h-\theta_{d-h}}
       {\theta_0-\theta_d}
 =
  \frac{\theta^*_h-\theta^*_{d-h}}
       {\theta^*_0-\theta^*_d}
   \qquad\qquad (0 \leq h \leq d).
\]
\end{lemma}

\begin{definition}         \label{def:varth}
Referring to Notation \ref{notation},
for $1 \leq i \leq d$ we have
\[
 \sum_{h=0}^{i-1}
  \frac{\theta_h-\theta_{d-h}}
       {\theta_0-\theta_d}
 =
 \sum_{h=0}^{i-1} \frac{\theta^*_h-\theta^*_{d-h}}
       {\theta^*_0-\theta^*_d}.
\]
We denote this common value by $\vartheta_i$.
We observe that $\vartheta_1=1$ and
$\vartheta_i=\vartheta_{d-i+1}$ for $1 \leq i \leq d$.
\end{definition}

\begin{lemma}         \label{eq:PA3}         \samepage
Referring to Notation \ref{notation} and
Definition \ref{def:varth},
the following hold for $1 \leq i \leq d$.
\begin{eqnarray}
 \varphi_i &=& \phi_1 \vartheta_i +
  (\theta^*_i-\theta^*_0)(\theta_{i-1}-\theta_d),    \label{eq:1}  \\
 \varphi_{d-i+1} &=&  \phi_1 \vartheta_i + 
  (\theta^*_{d-i+1}-\theta^*_0)(\theta_{d-i}-\theta_d),    \label{eq:2}  \\
 \varphi_i &=& \phi_d \vartheta_i +
  (\theta_i-\theta_0)(\theta^*_{i-1}-\theta^*_d),   \label{eq:3} \\
 \varphi_{d-i+1} &=& \phi_d \vartheta_i +
  (\theta_{d-i+1}-\theta_0)(\theta^*_{d-i}-\theta^*_d),   \label{eq:4}  \\
 \phi_i &=& \varphi_1 \vartheta_i +
  (\theta^*_i-\theta^*_0)(\theta_{d-i+1}-\theta_0),   \label{eq:5}  \\
 \phi_{d-i+1} &=& \varphi_1 \vartheta_i +
  (\theta^*_{d-i+1}-\theta^*_0)(\theta_i-\theta_0),   \label{eq:6}  \\
 \phi_i &=& \varphi_d \vartheta_i +
  (\theta_{d-i}-\theta_d)(\theta^*_{i-1}-\theta^*_d),  \label{eq:7}  \\
 \phi_{d-i+1} &=& \varphi_d \vartheta_i +
  (\theta_{i-1}-\theta_d)(\theta^*_{d-i}-\theta^*_d).  \label{eq:8}
\end{eqnarray}
\end{lemma}

\begin{proof}
Apply $D_4$ to the equation (PA3) from Theorem \ref{thm:classify}, and use
Lemma \ref{lem:D4}.
\end{proof}

\section{The relatives and affine transformations of a Leonard system}

\indent
Let $\Phi$ denote a Leonard system in $\cal A$.
In this section we give, for each relative of $\Phi$, 
necessary and sufficient conditions for it to be affine isomorphic to $\Phi$.
Recall that by Theorem \ref{thm:classify},
two Leonard systems are isomorphic if and only if
they have the same parameter array.

\begin{lemma}          \label{lem:actD4}            \samepage
Let $\Phi$ and $\Phi'$ denote Leonard systems over $\mathbb{K}$
which are affine isomorphic.
Then $\Phi^g$ and ${\Phi'}^g$ are affine isomorphic
for all $g \in D_4$.
\end{lemma}

\begin{proof}
Routine.
\end{proof}

\begin{proposition}           \label{prop:trivial}        \samepage
Referring to Notation \ref{notation}, let 
$\xi, \zeta, \xi^*, \zeta^*$ 
denote scalars in $\mathbb{K}$ with $\xi, \xi^*$ nonzero. 
Then $\Phi$ is isomorphic to the Leonard system (\ref{eq:affPhi})
if and only if $\xi=1$, $\zeta=0$, $\xi^*=1$, $\zeta^*=0$.
\end{proposition}

\begin{proof}
Suppose that $\Phi$ is isomorphic to the Leonard system (\ref{eq:affPhi}).
Then these Leonard systems have the same parameter array.
These parameter arrays are given in Notation \ref{notation} and 
(\ref{eq:affparam}); 
comparing them we find $\xi \theta_i+\zeta=\theta_i$ for $0 \leq i \leq d$.
Setting $i=0$, $i=1$ in this equation we find $\xi=1$, $\zeta=0$.
Similarly we find $\xi^*=1$, $\zeta^*=0$.
This proves the result in one direction
and the other direction is clear.
\end{proof}

\begin{lemma}          \label{lem:down}    \samepage
Referring to Notation \ref{notation},
let $\xi, \zeta, \xi^*, \zeta^*$ denote scalars in $\mathbb{K}$ with 
$\xi, \xi^*$ nonzero. 
Then $\Phi^\downarrow$ is isomorphic to the Leonard system (\ref{eq:affPhi})
if and only if
\begin{eqnarray}
 \theta_{i} &=& \xi \theta_i+\zeta 
       \qquad\qquad  (0 \leq i \leq d),   \label{eq:(i)aux1}  \\
 \theta^*_{d-i} &=& \xi^* \theta^*_i + \zeta^*
       \qquad\qquad (0 \leq i \leq d),    \label{eq:(i)aux2}  \\
  \phi_{d-i+1} &=& \xi\xi^* \varphi_i
       \qquad\qquad (1 \leq i \leq d),     \label{eq:(i)aux3}   \\
  \varphi_{d-i+1} &=& \xi\xi^* \phi_i 
       \qquad\qquad (1 \leq i \leq d).
           \label{eq:(i)aux4}
\end{eqnarray}
\end{lemma}

\begin{proof}
Compare the parameter array of $\Phi^{\downarrow}$ from Lemma \ref{lem:D4},
with the parameter array (\ref{eq:affparam}).
\end{proof}

\begin{proposition}           \label{prop:down}   \samepage
Referring to Notation \ref{notation},
the following (i)--(iii) are equivalent.
\begin{itemize}
\item[(i)]
$\Phi^{\downarrow}$ is affine isomorphic to $\Phi$.
\item[(ii)]
$\varphi_1=-\phi_d$ and $\varphi_d=-\phi_1$.
\item[(iii)]
$\varphi_i=-\phi_{d-i+1}$ for $1 \leq i \leq d$ and
$\theta^*_i+\theta^*_{d-i}$ is independent of $i$ for $0 \leq i \leq d$.
\end{itemize}
Suppose (i)--(iii) hold.
Then $\Phi^{\downarrow}$ is isomorphic to (\ref{eq:affPhi}) with
$\xi=1$, 
$\zeta=0$,
$\xi^*=-1$, and 
$\zeta^*$ equal to the common value of $\theta^*_i+\theta^*_{d-i}$.
\end{proposition}

\begin{proof}
(i)$\Rightarrow$(ii):
By Definition \ref{def:affiso} there exist scalars 
$\xi, \zeta, \xi^*, \zeta^*$ in $\mathbb{K}$ with $\xi, \xi^*$ nonzero such that 
$\Phi^\downarrow$ is isomorphic to the Leonard system (\ref{eq:affPhi}). 
Now (\ref{eq:(i)aux1})--(\ref{eq:(i)aux4}) hold by Lemma \ref{lem:down}.
Setting $i=0$, $i=1$ in (\ref{eq:(i)aux1}) we find $\xi=1$, $\zeta=0$.
Setting $i=0$, $i=d$ in (\ref{eq:(i)aux2}) we find $\xi^*=-1$.
Setting $i=1$, $i=d$ in (\ref{eq:(i)aux3}) and using $\xi=1$, $\xi^*=-1$ we find
$\varphi_1 = - \phi_d$ and $\varphi_d = - \phi_1$.

(ii)$\Rightarrow$(iii):
By (\ref{eq:1}), (\ref{eq:8}) and $\varphi_d=-\phi_1$,
\begin{equation}               \label{eq:vphi1phid1}
  \varphi_i+\phi_{d-i+1} 
 = (\theta_{i-1}-\theta_d)(\theta^*_i+\theta^*_{d-i}-\theta^*_0-\theta^*_d)
   \qquad\qquad (1 \leq i \leq d).
\end{equation}                  
By (\ref{eq:3}), (\ref{eq:6}) and $\varphi_1=-\phi_d$,
\begin{equation}                \label{eq:vphi1phid2}
  \varphi_i+\phi_{d-i+1} 
 = (\theta_{i}-\theta_0)(\theta^*_{i-1}+\theta^*_{d-i+1}-\theta^*_0-\theta^*_d)
   \qquad\qquad (1 \leq i \leq d).
\end{equation}
Replacing $i$ by $i+$1 in (\ref{eq:vphi1phid2}) and comparing the result
with (\ref{eq:vphi1phid1}) we find
\[
  \frac{\varphi_i + \phi_{d-i+1}}
       {\theta_{i-1} - \theta_d}
= \frac{\varphi_{i+1}+\phi_{d-i}}
       {\theta_{i+1}-\theta_0}
   \qquad\qquad  (1 \leq i  \leq d-1).
\]
From this and since $\varphi_1+\phi_d=0$ we find
$\varphi_i + \phi_{d-i+1}=0$ for $1 \leq i \leq d$.
Evaluating (\ref{eq:vphi1phid1}) using this we find
$\theta^*_i + \theta^*_{d-i}$ is independent of $i$ for $0 \leq i \leq d$.

(iii)$\Rightarrow$(i):
Let $\zeta^*$ denote the common value of $\theta^*_i+\theta^*_{d-i}$, and
let $\xi=1$, $\zeta=0$, $\xi^*=-1$.
Now (\ref{eq:(i)aux1})--(\ref{eq:(i)aux4}) hold so
$\Phi^\downarrow$ is isomorphic to (\ref{eq:affPhi}) by 
Lemma \ref{lem:down}.
Now $\Phi^\downarrow$ is affine isomorphic to $\Phi$ in
view of Definition \ref{def:affiso}.
\end{proof}

\begin{proposition}           \label{prop:Down}   \samepage
Referring to Notation \ref{notation},
the following (i)--(iii) are equivalent.
\begin{itemize}
\item[(i)]
$\Phi^{\Downarrow}$ is affine isomorphic to $\Phi$.
\item[(ii)]
$\varphi_1=-\phi_1$ and $\varphi_d=-\phi_d$.
\item[(iii)]
$\varphi_i=-\phi_{i}$ for $1 \leq i \leq d$ and
$\theta_i+\theta_{d-i}$ is independent of $i$ for $0 \leq i \leq d$.
\end{itemize}
Suppose (i)--(iii) hold.
Then $\Phi^{\Downarrow}$ is isomorphic to (\ref{eq:affPhi}) with
$\xi=-1$, 
$\zeta$ equal to the common value of $\theta_i+\theta_{d-i}$, 
$\xi^*=1$, and 
$\zeta^*=0$.
\end{proposition}

\begin{proof}
By Lemma \ref{lem:actD4} (with $g=*$) and since 
$\Downarrow\,$$*\,$$=\,$$*\,$$\downarrow$
we find
$\Phi^{\Downarrow}$ is affine 
isomorphic to $\Phi$ if and only if $\Phi^{*\downarrow}$ is affine 
isomorphic to $\Phi^*$.
Now apply Proposition \ref{prop:down} to $\Phi^*$ and use Lemma \ref{lem:D4}.
\end{proof}

\begin{lemma}          \label{lem:star}    \samepage
Referring to Notation \ref{notation},
let $\xi, \zeta, \xi^*, \zeta^*$ denote scalars in $\mathbb{K}$ with 
$\xi, \xi^*$ nonzero. 
Then $\Phi^*$ is isomorphic to the Leonard system (\ref{eq:affPhi})
if and only if
\begin{eqnarray}
  \theta^*_{i} &=& \xi \theta_i+\zeta    \qquad  (0 \leq i \leq d), 
            \label{eq:(iv)aux1}  \\
  \theta_{i} &=& \xi^* \theta^*_i + \zeta^*  \qquad (0 \leq i \leq d),
            \label{eq:(iv)aux2}  \\
  \varphi_{i} &=& \xi\xi^* \varphi_i    \qquad (1 \leq i \leq d), 
           \label{eq:(iv)aux3}   \\
  \phi_{d-i+1} &=& \xi\xi^* \phi_i   \qquad (1 \leq i \leq d).
           \label{eq:(iv)aux4}
\end{eqnarray}
\end{lemma}

\begin{proof}
Compare the parameter array of $\Phi^{*}$ from Lemma \ref{lem:D4},
with the parameter array (\ref{eq:affparam}).
\end{proof}

\begin{proposition}           \label{prop:star}   \samepage
Referring to Notation \ref{notation},
the following (i)--(iv) are equivalent.
\begin{itemize}
\item[(i)]
$\Phi^{*}$ is affine isomorphic to $\Phi$.
\item[(ii)] 
$\phi_1=\phi_d$.
\item[(iii)]
$\phi_i=\phi_{d-i+1}$ for $1 \leq i \leq d$.
\item[(iv)]
$(\theta^*_i-\theta^*_0)(\theta_i-\theta_0)^{-1}$ is independent of $i$
for $1 \leq i \leq d$.
\end{itemize}
Suppose (i)--(iv) hold.
Then $\Phi^{*}$ is isomorphic to (\ref{eq:affPhi}) with
$\xi$ equal to the common value of 
$(\theta^*_i-\theta^*_0)(\theta_i-\theta_0)^{-1}$,
$\zeta=\theta^*_0-\xi \theta_0$, 
$\xi^*=\xi^{-1}$, and 
$\zeta^*=\theta_0-\xi^* \theta^*_0$.
\end{proposition}

\begin{proof}
(i)$\Rightarrow$(ii):
By Definition \ref{def:affiso} there exist scalars 
$\xi, \zeta, \xi^*, \zeta^*$ in $\mathbb{K}$ with 
$\xi, \xi^*$ nonzero such that 
$\Phi^*$ is isomorphic to the Leonard system (\ref{eq:affPhi}). 
Now (\ref{eq:(iv)aux1})--(\ref{eq:(iv)aux4}) hold by Lemma \ref{lem:star}.
By (\ref{eq:(iv)aux3}) we find $\xi\xi^*=1$.
Setting $i=1$ in (\ref{eq:(iv)aux4}) and using $\xi\xi^*=1$ we find
$\phi_1=\phi_d$.

(ii)$\Rightarrow$(iv):
For $0 \leq i \leq d$ define
$\eta_i = (\theta^*_i-\theta^*_0)(\theta_d-\theta_0)- 
          (\theta_i-\theta_0)(\theta^*_d-\theta^*_0)$
and observe $\eta_0=0$. 
We show $\eta_i=0$ for $1 \leq i \leq d$.
By (\ref{eq:1}), (\ref{eq:3}) and since $\phi_1=\phi_d$,
\[
 (\theta^*_i-\theta^*_0)(\theta_{i-1}-\theta_d)
 =(\theta_i-\theta_0)(\theta^*_{i-1}-\theta^*_d)
  \qquad\qquad     (1 \leq i \leq d).
\]
In this equation we rearrange terms to get
\[
 \eta_i (\theta_{i-1}-\theta_d)= \eta_{i-1}(\theta_i-\theta_0) 
  \qquad\qquad      (1 \leq i \leq d).
\]
By this and since $\eta_0 =0$ we find $\eta_i=0$ for $1 \leq i \leq d$.
The result follows.

(iv)$\Rightarrow$(iii):
Let $i$ be given. Since
$(\theta^*_i-\theta^*_0)(\theta_i-\theta_0)^{-1}$ is independent of $i$,
\[
  (\theta^*_i-\theta^*_0)(\theta_{d-i+1}-\theta_0)
 = (\theta^*_{d-i+1}-\theta^*_0)(\theta_i-\theta_0).
\]
Comparing (\ref{eq:5}) and (\ref{eq:6}) using this we find 
$\phi_i = \phi_{d-i+1}$.

(iii)$\Rightarrow$(ii):
Clear.

(iii),$\,$(iv)$\Rightarrow$(i):
Let $\xi$ denote the common value of 
$(\theta^*_i-\theta^*_0)(\theta_i-\theta_0)^{-1}$ and set
$\xi^*=\xi^{-1}$,
$\zeta=\theta^*_0 - \xi\theta_0$,
$\zeta^*=\theta_0-\xi^*\theta^*_0$.
Then (\ref{eq:(iv)aux1})--(\ref{eq:(iv)aux4}) hold so
$\Phi^*$ is isormorphic to (\ref{eq:affPhi}) by Lemma \ref{lem:star}.
Now $\Phi^*$ is affine isomorphic to $\Phi$ in
view of Definition \ref{def:affiso}.
\end{proof}

\begin{proposition}           \label{prop:downDownstar}    \samepage
Referring to Notation \ref{notation},
the following (i)--(iv) are equivalent.
\begin{itemize}
\item[(i)]
$\Phi^{\downarrow\Downarrow*}$ is affine isomorphic to $\Phi$.
\item[(ii)]
$\varphi_1=\varphi_d$.
\item[(iii)]
$\varphi_i=\varphi_{d-i+1}$ for $1 \leq i \leq d$.
\item[(iv)]
$(\theta^*_{d-i}-\theta^*_d)(\theta_i-\theta_0)^{-1}$
is independent of $i$ for $1 \leq i \leq d$.
\end{itemize}
Suppose (i)--(iv) hold.
Then $\Phi^{\downarrow\Downarrow*}$ is isomorphic to (\ref{eq:affPhi}) with
$\xi$ equal to the common value of 
$(\theta^*_{d-i}-\theta^*_d)(\theta_i-\theta_0)^{-1}$,
$\zeta=\theta^*_d-\xi \theta_0$,
$\xi^*=\xi^{-1}$, and 
$\zeta^*=\theta_0-\xi^* \theta^*_d$.
\end{proposition}

\begin{proof}
By Lemma \ref{lem:actD4} (with $g = \,\,\downarrow$) and since
$\Downarrow\,$$*\,$$\downarrow\,$$=\,$$*$
we find that
$\Phi^{\downarrow\Downarrow*}$ is affine isomorphic to $\Phi$
if and only if $\Phi^{\downarrow*}$ is affine isomorphic to 
$\Phi^{\downarrow}$. 
Now apply Proposition \ref{prop:star} to
$\Phi^{\downarrow}$ and use Lemma \ref{lem:D4}.
\end{proof}

\begin{lemma}          \label{lem:downDown}    \samepage
Referring to Notation \ref{notation},
let $\xi, \zeta, \xi^*, \zeta^*$ denote scalars in $\mathbb{K}$ with 
$\xi, \xi^*$ nonzero. 
Then $\Phi^{\downarrow\Downarrow}$ is isomorphic to the Leonard system 
(\ref{eq:affPhi}) if and only if
\begin{eqnarray}
  \theta_{d-i} &=& \xi \theta_i+\zeta 
       \qquad\qquad  (0 \leq i \leq d),       \label{eq:(iii)aux1}  \\
  \theta^*_{d-i} &=& \xi^* \theta^*_i + \zeta^* 
       \qquad\qquad (0 \leq i \leq d),        \label{eq:(iii)aux2}  \\
  \varphi_{d-i+1} &=& \xi\xi^* \varphi_i 
       \qquad\qquad (1 \leq i \leq d),        \label{eq:(iii)aux3}   \\
  \phi_{d-i+1} &=& \xi\xi^* \phi_i
       \qquad\qquad (1 \leq i \leq d).        \label{eq:(iii)aux4}
\end{eqnarray}
\end{lemma}

\begin{proof}
Compare the parameter array of $\Phi^{\downarrow\Downarrow}$ from 
Lemma \ref{lem:D4},
with the parameter array (\ref{eq:affparam}).
\end{proof}

\begin{proposition}           \label{prop:downDown}   \samepage
Referring to Notation \ref{notation},
the following (i)--(iv) are equivalent.
\begin{itemize}
\item[(i)]
$\Phi^{\downarrow\Downarrow}$ is affine isomorphic to $\Phi$.
\item[(ii)]
$\varphi_1=\varphi_d$ and $\phi_1=\phi_d$.
\item[(iii)]
$\varphi_i=\varphi_{d-i+1}$ and $\phi_i=\phi_{d-i+1}$ for $1 \leq i \leq d$.
\item[(iv)]
Each of $(\theta^*_i-\theta^*_0)(\theta_i-\theta_0)^{-1}$,
$(\theta^*_{d-i}-\theta^*_d)(\theta_i-\theta_0)^{-1}$ is independent of $i$
for $1 \leq i \leq d$.
\end{itemize}
Suppose (i)--(iv) hold.
Then each of $\theta_i+\theta_{d-i}$, $\theta^*_i+\theta^*_{d-i}$ is
independent of $i$ for $0 \leq i \leq d$.
Moreover
$\Phi^{\downarrow\Downarrow}$ is isomorphic to (\ref{eq:affPhi}) with
$\xi=-1$, 
$\xi^*=-1$, 
and $\zeta$ $($resp. $\zeta^*)$ equal to the common value of 
$\theta_i+\theta_{d-i}$ $($resp. $\theta^*_i+\theta^*_{d-i})$.
\end{proposition}

\begin{proof}
(i)$\Rightarrow$(ii):
By Definition \ref{def:affiso} there exist scalars 
$\xi, \zeta, \xi^*, \zeta^*$ in $\mathbb{K}$ with $\xi, \xi^*$ nonzero such that 
$\Phi^{\downarrow\Downarrow}$ is isomorphic to the Leonard system 
(\ref{eq:affPhi}). 
Now (\ref{eq:(iii)aux1})--(\ref{eq:(iii)aux4}) hold by Lemma \ref{lem:downDown}.
Setting $i=0$, $i=d$ in (\ref{eq:(iii)aux1}) we find $\xi=-1$.
Setting $i=0$, $i=d$ in (\ref{eq:(iii)aux2}) we find $\xi^*=-1$.
Setting $i=1$ in (\ref{eq:(iii)aux3}), (\ref{eq:(iii)aux4})
and using $\xi=-1$, $\xi^*=-1$ we find
$\varphi_1 = \varphi_d$ and $\phi_1 = \phi_d$.

(ii)$\Leftrightarrow$(iii)$\Leftrightarrow$(iv):
Follows from Propositions \ref{prop:star} and \ref{prop:downDownstar}.

(iii),$\,$(iv)$\Rightarrow$(i):
We first show that $\theta^*_i+\theta^*_{d-i}$
is independent of $i$ for $0 \leq i \leq d$.
By assumption
\begin{equation}             \label{eq:thsiths01}
   \frac{\theta^*_i-\theta^*_0}
        {\theta_i-\theta_0}
 = \frac{\theta^*_d-\theta^*_0}
        {\theta_d-\theta_0}
   \qquad\qquad (1 \leq i \leq d),
\end{equation}
and 
\begin{equation}           \label{eq:thsiths02}
   \frac{\theta^*_{d-i}-\theta^*_d}
        {\theta_i-\theta_0}
 = \frac{\theta^*_0-\theta^*_d}
        {\theta_d-\theta_0}
   \qquad\qquad (1 \leq i \leq d).
\end{equation}
Adding (\ref{eq:thsiths01}), (\ref{eq:thsiths02})  we find
$\theta^*_i+\theta^*_{d-i}=\theta^*_0+\theta^*_d$
for $1 \leq i \leq d$.
Therefore
$\theta^*_i+\theta^*_{d-i}$ is independent of $i$ for $0 \leq i \leq d$.
Next we show that $\theta_i+\theta_{d-i}$ is independent of $i$
for $0 \leq i \leq d$.
Rearranging the terms in (\ref{eq:thsiths01}) we find that for
$1 \leq i \leq d$,
\[
  \frac{\theta_i+\theta_{d-i}-\theta_0-\theta_d}
       {\theta_0-\theta_d} 
= \frac{\theta^*_i+\theta^*_{d-i}-\theta^*_0-\theta^*_d}
       {\theta^*_0-\theta^*_d}.
\]
In the above equation the numerator on the right is zero so
the numerator on the left is zero. Therefore
$\theta_i+\theta_{d-i}$ is independent of $i$ for $0 \leq i \leq d$.
Now let $\zeta$ (resp. $\zeta^*$) denote the common value of 
$\theta_i+\theta_{d-i}$ (resp. $\theta^*_i+\theta^*_{d-i}$), and
let $\xi=-1$, $\xi^*=-1$.
Then (\ref{eq:(iii)aux1})--(\ref{eq:(iii)aux4}) hold so
$\Phi^{\downarrow\Downarrow}$ is isomorpic to (\ref{eq:affPhi})
by Lemma \ref{lem:downDown}.
Now $\Phi^{\downarrow\Downarrow}$ is affine isomorphic to $\Phi$ in
view of Definition \ref{def:affiso}.
\end{proof}

\begin{lemma}          \label{lem:downstar}    \samepage
Referring to Notation \ref{notation},
let $\xi, \zeta, \xi^*, \zeta^*$ denote scalars in $\mathbb{K}$ with 
$\xi, \xi^*$ nonzero. 
Then $\Phi^{\downarrow*}$ is isomorphic to the Leonard system (\ref{eq:affPhi})
if and only if
\begin{eqnarray}
  \theta^*_{d-i} &=& \xi \theta_i+\zeta    \qquad  (0 \leq i \leq d), 
            \label{eq:(v)aux1}  \\
  \theta_{i} &=& \xi^* \theta^*_i + \zeta^*   \qquad (0 \leq i \leq d),
            \label{eq:(v)aux2}  \\
  \phi_{d-i+1} &=& \xi\xi^* \varphi_i    \qquad (1 \leq i \leq d), 
           \label{eq:(v)aux3}   \\
  \varphi_{i} &=&  \xi\xi^* \phi_i     \qquad (1 \leq i \leq d).
           \label{eq:(v)aux4}
\end{eqnarray}
\end{lemma}

\begin{proof}
Compare the parameter array of $\Phi^{\downarrow*}$ from Lemma \ref{lem:D4},
with the parameter array (\ref{eq:affparam}).
\end{proof}

\begin{proposition}           \label{prop:downstar}    \samepage
Referring to Notation \ref{notation},
the following (i)--(iii) are equivalent.
\begin{itemize}
\item[(i)]
$\Phi^{\downarrow*}$ is affine isomorphic to $\Phi$.
\item[(ii)]
$\varphi_1=\varphi_d=-\phi_1=-\phi_d$.
\item[(iii)]
$\varphi_i$, $\varphi_{d-i+1}$, $-\phi_i$, $-\phi_{d-i+1}$ 
coincide for $1 \leq i \leq d$
and each of 
$(\theta^*_{i}-\theta^*_0)(\theta_i-\theta_0)^{-1}$,
$(\theta^*_{d-i}-\theta^*_d)(\theta_i-\theta_0)^{-1}$
is independent of $i$ for $1 \leq i \leq d$.
\end{itemize}
Suppose (i)--(iii) hold.
Then $\Phi^{\downarrow*}$ is isomorphic to (\ref{eq:affPhi}) with
$\xi$ equal to the common value of 
$(\theta^*_{d-i}-\theta^*_d)(\theta_i-\theta_0)^{-1}$,
$\zeta=\theta^*_d-\xi \theta_0$, 
$\xi^*=- \xi^{-1}$, and 
$\zeta^*=\theta_0-\xi^* \theta^*_0$.
\end{proposition}

\begin{proof}
(i)$\Rightarrow$(ii):
By Definition \ref{def:affiso} there exist scalars 
$\xi, \zeta, \xi^*, \zeta^*$ in $\mathbb{K}$ with 
$\xi, \xi^*$ nonzero such that 
$\Phi^{\downarrow*}$ is isomorphic to the Leonard system (\ref{eq:affPhi}). 
Now (\ref{eq:(v)aux1})--(\ref{eq:(v)aux4}) hold by Lemma \ref{lem:downstar}.
Setting $i=0$, $i=d$ in (\ref{eq:(v)aux1}) we find 
$\xi(\theta_0-\theta_d)=\theta^*_d-\theta^*_0$.
Setting $i=0$, $i=d$ in (\ref{eq:(v)aux2}) we find 
$\xi^*(\theta^*_0-\theta^*_d)=\theta_0-\theta_d$.
By these comments $\xi\xi^*=-1$.
Setting $i=1$, $i=d$ in (\ref{eq:(v)aux3}) and using $\xi\xi^*=-1$
we find $\varphi_1=-\phi_d$ and $\varphi_d=-\phi_1$.
Setting $i=1$ in (\ref{eq:(v)aux4}) and using $\xi\xi^*=-1$
we find $\varphi_1=-\phi_1$.

(ii)$\Leftrightarrow$(iii):
Follows from Propositions \ref{prop:down}, \ref{prop:Down} and
\ref{prop:downDown}.

(ii),$\,$(iii)$\Rightarrow$(i):
Let $\xi$ denote the common value of 
$(\theta^*_{d-i} - \theta^*_d)(\theta_i-\theta_0)^{-1}$, and let
$\xi^* = -\xi^{-1}$, $\zeta=\theta^*_d-\xi \theta_0$,
$\zeta^*=\theta_0-\xi^*\theta^*_0$.
Then (\ref{eq:(v)aux1})--(\ref{eq:(v)aux4}) hold so
$\Phi^{\downarrow*}$ is isormorphic to (\ref{eq:affPhi})
by Lemma \ref{lem:downstar}.
Now $\Phi^{\downarrow*}$ is affine isomorphic to $\Phi$ in
view of Definition \ref{def:affiso}.
\end{proof}

\begin{proposition}           \label{prop:Downstar}    \samepage
Referring to Notation \ref{notation},
the following (i)--(iii) are equivalent.
\begin{itemize}
\item[(i)]
$\Phi^{\Downarrow*}$ is affine isomorphic to $\Phi$.
\item[(ii)]
$\varphi_1=\varphi_d=-\phi_1=-\phi_d$.
\item[(iii)]
$\varphi_i$, $\varphi_{d-i+1}$, $-\phi_i$, $-\phi_{d-i+1}$ 
coincide for $1 \leq i \leq d$ and each of
$(\theta^*_{i}-\theta^*_0)(\theta_i-\theta_0)^{-1}$,
$(\theta^*_{d-i}-\theta^*_d)(\theta_i-\theta_0)^{-1}$
is independent of $i$ for $1 \leq i \leq d$.
\end{itemize}
Suppose (i)--(iii) hold.
Then $\Phi^{\Downarrow*}$ is affine isomorphic to (\ref{eq:affPhi}) with
$\xi$ equal to the common value of 
$(\theta^*_{i}-\theta^*_0)(\theta_i-\theta_0)^{-1}$,
$\zeta=\theta^*_0-\xi \theta_0$, $\xi^*=- \xi^{-1}$, and 
$\zeta^*=\theta_0-\xi^* \theta^*_d$.
\end{proposition}

\begin{proof}
By Lemma \ref{lem:actD4} (with $g = \,\downarrow$)
and since 
$\Downarrow\,$$*\,$$\downarrow\,$$=\,$$*\,$$=\,$$\downarrow\,$$\downarrow\,$$*$
we find that
$\Phi^{\Downarrow*}$ is affine isomorphic to $\Phi$
if and only if
$(\Phi^\downarrow)^{\downarrow*}$ is affine isomorphic to $\Phi^{\downarrow}$.
Now apply Proposition \ref{prop:downstar} to
$\Phi^{\downarrow}$ and use Lemma \ref{lem:D4}.
\end{proof}

\section{Proof of Theorem \ref{thm:main}}

\indent
In this section we prove Theorem \ref{thm:main}.

\begin{proofof}{Theorem \ref{thm:main}}
(i): 
Observe that $\Phi^g$ is isomorphic to $\Phi$ for all $g \in D_4$ by
Propositions 
\ref{prop:down}, \ref{prop:Down}, \ref{prop:star}, \ref{prop:downDownstar},
\ref{prop:downDown}, \ref{prop:downstar} and \ref{prop:Downstar}.

(ii):
Since $\varphi_1=\varphi_d$ and $\phi_1=\phi_d$,
the Leonard systems 
$\Phi^{*}$,
$\Phi^{\downarrow\Downarrow*}$,
$\Phi^{\downarrow\Downarrow}$
are affine isomorphic to $\Phi$ by Propositions  
\ref{prop:star}, 
\ref{prop:downDownstar},
\ref{prop:downDown}
respectively.
Therefore
$\Phi$,
$\Phi^{*}$,
$\Phi^{\downarrow\Downarrow*}$,
$\Phi^{\downarrow\Downarrow}$
are contained in a common affine isomorphism class.
By this and Lemma \ref{lem:actD4} the Leonard systems
$\Phi^{\downarrow}$, $\Phi^{\Downarrow}$, 
$\Phi^{\downarrow*}$, $\Phi^{\Downarrow*}$ are contained in
a common isomorphism class.
The above affine isomorphism classes are distinct;
indeed $\Phi^\downarrow$ is not affine isomorphic to $\Phi$
by Proposition  \ref{prop:down}  and since
$\varphi_1 \neq -\phi_d$. The result follows.

(iii):
By Propostion \ref{prop:downDownstar} the Leonard system $\Phi$ 
is affine isomorphic to $\Phi^{\downarrow\Downarrow*}$.
By Propositions
\ref{prop:down},
\ref{prop:Down},
\ref{prop:downDown},
\ref{prop:star},
\ref{prop:downstar},
\ref{prop:Downstar},
$\Phi$ is not affine isomorphic to any of
$\Phi^{\downarrow}$,
$\Phi^{\Downarrow}$,
$\Phi^{\downarrow\Downarrow}$,
$\Phi^{*}$,
$\Phi^{\downarrow*}$,
$\Phi^{\Downarrow*}$.
The result follows from these comments in view of Lemma \ref{lem:actD4}
and (\ref{eq:relation1}), (\ref{eq:relation2}).

(iv):
By Propostion \ref{prop:star} the Leonard system $\Phi$ 
is affine isomorphic to $\Phi^*$.
By Propositions
\ref{prop:down},
\ref{prop:Down},
\ref{prop:downDown},
\ref{prop:downstar},
\ref{prop:Downstar},
\ref{prop:downDownstar},
$\Phi$ is not affine isomorphic to any of
$\Phi^{\downarrow}$,
$\Phi^{\Downarrow}$,
$\Phi^{\downarrow\Downarrow}$,
$\Phi^{\downarrow*}$,
$\Phi^{\Downarrow*}$,
$\Phi^{\downarrow\Downarrow*}$.
The result follows from these comments in view of Lemma \ref{lem:actD4}
and (\ref{eq:relation1}), (\ref{eq:relation2}).

(v):
By Propostion \ref{prop:Down} the Leonard system $\Phi$ 
is affine isomorphic to $\Phi^{\Downarrow}$.
By Propositions
\ref{prop:down},
\ref{prop:downDown},
\ref{prop:star},
\ref{prop:downstar},
\ref{prop:Downstar},
\ref{prop:downDownstar},
$\Phi$ is not affine isomorphic to any of
$\Phi^{\downarrow}$,
$\Phi^{\downarrow\Downarrow}$,
$\Phi^{*}$,
$\Phi^{\downarrow*}$,
$\Phi^{\Downarrow*}$,
$\Phi^{\downarrow\Downarrow*}$.
The result follows from these comments in view of Lemma \ref{lem:actD4}
and (\ref{eq:relation1}), (\ref{eq:relation2}).

(vi):
By Propostion \ref{prop:down} the Leonard system $\Phi$ 
is affine isomorphic to $\Phi^{\downarrow}$.
By Propositions
\ref{prop:Down},
\ref{prop:downDown},
\ref{prop:star},
\ref{prop:downstar},
\ref{prop:Downstar},
\ref{prop:downDownstar},
$\Phi$ is not affine isomorphic to any of
$\Phi^{\Downarrow}$,
$\Phi^{\downarrow\Downarrow}$,
$\Phi^{*}$,
$\Phi^{\downarrow*}$,
$\Phi^{\Downarrow*}$,
$\Phi^{\downarrow\Downarrow*}$.
The result follows from these comments in view of Lemma \ref{lem:actD4}
and (\ref{eq:relation1}), (\ref{eq:relation2}).

(vii):
By Propositions
\ref{prop:down},
\ref{prop:Down},
\ref{prop:star},
\ref{prop:downDown},
\ref{prop:downstar},
\ref{prop:Downstar},
\ref{prop:downDownstar},
$\Phi$ is not affine isomorphic to any of
$\Phi^{\downarrow}$,
$\Phi^{\Downarrow}$,
$\Phi^{*}$,
$\Phi^{\downarrow\Downarrow}$,
$\Phi^{\downarrow*}$,
$\Phi^{\Downarrow*}$,
$\Phi^{\downarrow\Downarrow*}$.
The result follows from this and  Lemma \ref{lem:actD4}.
\end{proofof}

\section{The parameters $a_i$ and $a^*_i$}

\indent
It turns out that for some of the cases of Theorem \ref{thm:main}
there is a natural interpretation in terms of the parameters
$a_i$ and $a^*_i$ \cite[Definition 2.5]{T:Leonard}.
In this section we explain the situation. 
We start with a definition.

\begin{definition}  \cite[Definition 2.5]{T:Leonard}  \label{def:ai}
Referring to Notation \ref{notation},
for $0 \leq i \leq d$ we define scalars
\[
  a_i:=\text{tr}(E^*_iA), \qquad\qquad
  a^*_i:=\text{tr}(E_iA^*). 
\]
\end{definition}

\begin{lemma}            \label{lem:bip}       \samepage
Referring to Notation \ref{notation} and Definition \ref{def:ai},
the following (i), (ii) are equaivalent.
\begin{itemize}
\item[(i)]
$a_i$ is independent of $i$ for $0 \leq i \leq d$.
\item[(ii)]
The equivalent conditions (i)--(iii) hold in Proposition \ref{prop:Down}.
\end{itemize}
Suppose (i), (ii) hold. 
Then the common value of $\theta_i+\theta_{d-i}$ is twice the
common value of $a_i$.
\end{lemma}

\begin{proof}
Follows from \cite[Theorem 1.5]{NT:balanced} and Proposition \ref{prop:Down}.
\end{proof}

\begin{lemma}            \label{lem:dualbip}     \samepage
Referring to Notation \ref{notation} and Definition \ref{def:ai},
the following (i), (ii) are equaivalent.
\begin{itemize}
\item[(i)]
$a^*_i$ is independent of $i$ for $0 \leq i \leq d$.
\item[(ii)]
The equivalent conditions (i)--(iii) hold in Proposition \ref{prop:down}.
\end{itemize}
Suppose (i), (ii) hold. 
Then the common value of $\theta^*_i+\theta^*_{d-i}$ is twice the
common value of $a^*_i$.
\end{lemma}

\begin{proof}
Follows from \cite[Theorem 1.6]{NT:balanced} and Proposition \ref{prop:down}.
\end{proof}

\begin{theorem}                               \samepage
Referring to Notation \ref{notation} and Definition \ref{def:ai},
the following (i)--(iv) hold.
\begin{itemize}
\item[(i)] 
In Case (i) of Theorem \ref{thm:main},
each of $a_i$, $a^*_i$ is independent of $i$ for
$0 \leq i \leq d$.
\item[(ii)]
In Case (v) of Theorem \ref{thm:main},
 $a_i$ is independent of $i$ for $0 \leq i \leq d$
but $a^*_i$ is not independent of $i$ for $0 \leq i \leq d$.
\item[(iii)]
In Case (vi) of Theorem \ref{thm:main},
$a^*_i$ is independent of $i$ for $0 \leq i \leq d$
but $a_i$ is not independent of $i$ for $0 \leq i \leq d$.
\item[(iv)]
In the remaining cases of Theorem \ref{thm:main},
neither of $a_i$, $a^*_i$ is independent of $i$ for $0 \leq i \leq d$.
\end{itemize}
\end{theorem}

\begin{proof}
Follows from Theorem \ref{thm:main} and 
Lemmas \ref{lem:bip}, \ref{lem:dualbip}.
\end{proof}

\section{Affine transformations of a Leonard pair}

\indent
Let $A,A^*$ denote a Leonard pair in $\cal A$
and let $\xi, \zeta, \xi^*, \zeta^*$ denote scalars
in $\mathbb{K}$ with $\xi, \xi^*$ nonzero.
By Lemma \ref{lem:affine} and our comments below
Definition \ref{def:LS} the pair
\begin{equation}         \label{eq:affLP}
     \xi A + \zeta I,\;\; \xi^* A^* +\zeta^* I
\end{equation}
is a Leonard pair in $\cal A$.
We call (\ref{eq:affLP}) the {\em affine transformation} of $A,A^*$
associated with $\xi, \zeta, \xi^*, \zeta^*$.
In this section we find necessary and sufficient conditions for
the Leonard pair (\ref{eq:affLP}) to be isomorphic to $A,A^*$.
We also find necessary and sufficient conditions for
the Leonard pair (\ref{eq:affLP}) to be isomorphic to 
the Leonard pair $A^*,A$.

\begin{notation}        \label{notation2}
Let $A,A^*$ denote a Leonard pair in $\cal A$.
Let $\Phi=(A; \{E_i\}_{i=0}^d; A^*; \{E^*_i\}_{i=0}^d)$
denote a Leonard system associated with $A,A^*$ and let
 $(\{\theta_i\}_{i=0}^d; \{\theta^*_i\}_{i=0}^d;
   \{\varphi_i\}_{i=1}^d; \{\phi_i\}_{i=1}^d)$
denote the parameter array of $\Phi$.
To avoid trivialities we assume $d \geq 1$.
\end{notation}

\begin{proposition}          \label{prop:LPaff}    \samepage
Referring to Notation \ref{notation2},
the Leonard pair $A,A^*$ is isomorphic to the Leonard pair
(\ref{eq:affLP}) if and only if
at least one of  (i)--(iv) holds below.
\begin{itemize}
\item[(i)]
$\xi=1$, $\zeta=0$, $\xi^*=1$, $\zeta^*=0$.
\item[(ii)]
$\varphi_1=-\phi_d$, $\varphi_d=-\phi_1$, 
$\xi=1$, $\zeta=0$, $\xi^*=-1$, $\zeta^*=\theta^*_0+\theta^*_d$.
\item[(iii)]
$\varphi_1=-\phi_1$, $\varphi_d=-\phi_d$,
$\xi=-1$, $\zeta=\theta_0+\theta_d$, $\xi^*=1$, $\zeta^*=0$.
\item[(iv)]
$\varphi_1=\varphi_d$, $\phi_1=\phi_d$,
$\xi=-1$, $\zeta=\theta_0+\theta_d$, $\xi^*=-1$,
$\zeta^*=\theta^*_0+\theta^*_d$.
\end{itemize}
In this case precisely one of (i)--(iv) holds.
\end{proposition}

\begin{proof}
By \cite[Lemma 5.4]{T:canform} the Leonard systems associated 
with $A,A^*$ are
$\Phi$, $\Phi^\downarrow$, $\Phi^\Downarrow$, $\Phi^{\downarrow\Downarrow}$.
Therefore the Leonard pair $A,A^*$ is isomorphic to the Leonard pair
(\ref{eq:affLP}) if and only if at least one of
$\Phi$, $\Phi^\downarrow$, $\Phi^\Downarrow$, $\Phi^{\downarrow\Downarrow}$
is isomorphic to the Leonard system (\ref{eq:affPhi}).
By this and Propositions
\ref{prop:trivial}, \ref{prop:down}, \ref{prop:Down}, \ref{prop:downDown} 
we find
$A,A^*$ is isomorphic to (\ref{eq:affLP}) if and only if 
at least one of (i)--(iv) holds.
Assume that at least one of (i)--(iv) holds.
We show that precisely one of (i)--(iv) holds.
By way of contradiction assume that at least two of (i)--(iv) hold.
Then at least one of $2\xi$, $2\xi^*$ is zero,
forcing  $\text{Char}(\mathbb{K})=2$,
and at least one of
$\theta_0+\theta_d$, $\theta^*_0+\theta^*_d$ is zero,
forcing  $\text{Char}(\mathbb{K}) \neq 2$ and giving
a contradiction.
Therefore  precisely one of (i)--(iv) holds.
\end{proof}

\begin{proposition}          \label{prop:LPaffs}    \samepage
Referring to Notation \ref{notation2},
the Leonard pair $A^*,A$ is isomorphic to the Leonard pair
(\ref{eq:affLP}) if and only if
at least one of  (i)--(iv) holds below.
\begin{itemize}
\item[(i)]
$\phi_1=\phi_d$,
$\xi=(\theta^*_d-\theta^*_0)(\theta_d-\theta_0)^{-1}$, 
$\zeta=\theta^*_0-\xi\theta_0$, 
$\xi^*=\xi^{-1}$,
$\zeta^*=\theta_0-\xi^*\theta^*_0$.
\item[(ii)]
$\varphi_1=\varphi_d=-\phi_1=-\phi_d$,
$\xi=(\theta^*_0-\theta^*_d)(\theta_d-\theta_0)^{-1}$, 
$\zeta=\theta^*_d-\xi\theta_0$, 
$\xi^*=-\xi^{-1}$,
$\zeta^*=\theta_0-\xi^*\theta^*_0$.
\item[(iii)]
$\varphi_1=\varphi_d=-\phi_1=-\phi_d$,
$\xi=(\theta^*_d-\theta^*_0)(\theta_d-\theta_0)^{-1}$, 
$\zeta=\theta^*_0-\xi\theta_0$, 
$\xi^*=-\xi^{-1}$,
$\zeta^*=\theta_0-\xi^*\theta^*_d$.
\item[(iv)]
$\varphi_1=\varphi_d$,
$\xi=(\theta^*_0-\theta^*_d)(\theta_d-\theta_0)^{-1}$, 
$\zeta=\theta^*_d-\xi\theta_0$, 
$\xi^*=\xi^{-1}$,
$\zeta^*=\theta_0-\xi^*\theta^*_d$.
\end{itemize}
In this case precisely one of (i)--(iv) holds.
\end{proposition}

\begin{proof}
By \cite[Lemma 5.4]{T:canform} the Leonard systems associated 
with $A^*,A$ are 
$\Phi^*$, $\Phi^{\downarrow*}$, 
$\Phi^{\Downarrow*}$, $\Phi^{\downarrow\Downarrow*}$.
Therefore the Leonard pair $A^*,A$ is isomorphic to the Leonard pair
(\ref{eq:affLP}) if and only if at least one of
$\Phi^*$, $\Phi^{\downarrow*}$, 
$\Phi^{\Downarrow*}$, $\Phi^{\downarrow\Downarrow*}$
is isomorphic to the Leonard system (\ref{eq:affPhi}).
By this and Propositions
\ref{prop:star}, \ref{prop:downDownstar},
\ref{prop:downstar}, \ref{prop:Downstar}
we find
$A^*,A$ is isomorphic to (\ref{eq:affLP}) 
if and only if at least one of (i)--(iv) holds.
Assume that at least one of (i)--(iv) holds.
We show that precisely one of (i)--(iv) holds.
By way of contradiction assume that at least two of (i)--(iv) hold.
Then at least one of $\theta_0=\theta_d$, $\theta^*_0=\theta^*_d$ holds,
for a contradiction.
Therefore precisely one of (i)--(iv) holds.
\end{proof}

\medskip

The following is the main result of the paper.

\medskip

\begin{theorem}       \label{thm:LPmain}  
Referring to Notation \ref{notation2}, we set
$\alpha=(\theta^*_d-\theta^*_0)(\theta_d-\theta_0)^{-1}$.
\begin{itemize}
\item[(i)]
Assume $\varphi_1=\varphi_d=-\phi_1=-\phi_d$. 
Then $A,A^*$ is isomorphic to (\ref{eq:affLP}) if and only if
the sequence $\xi, \zeta, \xi^*, \zeta^*$ is listed in the following table.
\[
 \begin{array}{cccc}
  \;\;\; \xi \;\;\;  & \zeta  & 
    \;\;\; \xi^* \;\;\; & \zeta^* \\
 \hline
   1 & 0 & 1 &0 \\
   1 & 0 & -1 & \theta^*_0+\theta^*_d \\
  -1 & \theta_0+\theta_d & 1 & 0 \\
  -1 & \theta_0+\theta_d & -1 & \theta^*_0+\theta^*_d 
 \end{array}
\]
Moreover $A^*,A$ is isomorphic to (\ref{eq:affLP}) if and only if
the sequence $\xi, \zeta, \xi^*, \zeta^*$ is listed in the following table.
\[
 \begin{array}{cccc}
  \;\;\; \xi \;\;\; & \zeta & \;\;\; \xi^* \;\;\; & \zeta^* \\
 \hline
   \alpha & \theta^*_0 - \alpha \theta_0 & 
   \alpha^{-1} & \theta_0 - \alpha^{-1} \theta^*_0
\\
   - \alpha & \theta^*_d + \alpha \theta_0 &
     \alpha^{-1} & \theta_0 - \alpha^{-1} \theta^*_0
\\
   \alpha & \theta^*_0 - \alpha \theta_0 &
      - \alpha^{-1} & \theta_0 + \alpha^{-1} \theta^*_d
\\
  - \alpha & \theta^*_d + \alpha \theta_0 &
  - \alpha^{-1} & \theta_0 + \alpha^{-1} \theta^*_d
 \end{array}
\]
\item[(ii)]
Assume $\varphi_1=\varphi_d$, $\phi_1=\phi_d$ and $\varphi_1 \neq -\phi_1$.
Then $A,A^*$ is isomorphic to (\ref{eq:affLP}) if and only if
the sequence $\xi, \zeta, \xi^*, \zeta^*$  is listed in the following table.
\[
 \begin{array}{cccc}
  \;\;\; \xi \;\;\; & \zeta & \;\;\; \xi^* \;\;\; & \zeta^* \\
 \hline
   1 & 0 & 1 & 0 \\
  -1 & \theta_0+\theta_d & -1 & \theta^*_0+\theta^*_d 
 \end{array}
\]
Moreover $A^*,A$ is isomorphic to (\ref{eq:affLP}) if and only if
the sequence $\xi$, $\zeta$, $\xi^*$, $\zeta^*$  is listed in the following table.
\[
 \begin{array}{cccc}
   \;\;\; \xi \;\;\; & \zeta & \;\;\; \xi^* \;\;\; & \zeta^* \\
 \hline
   \alpha & \theta^*_0 - \alpha \theta_0 & 
   \alpha^{-1} & \theta_0 - \alpha^{-1} \theta^*_0
\\
  - \alpha & \theta^*_d + \alpha \theta_0 &
  - \alpha^{-1} & \theta_0 + \alpha^{-1} \theta^*_d
 \end{array}
\]
\item[(iii)] 
Assume $\varphi_1=\varphi_d$ and $\phi_1 \neq \phi_d$.
Then $A,A^*$ is isomorphic to (\ref{eq:affLP}) if and only if
$\xi=1$, $\zeta=0$, $\xi^*=1$, $\zeta^*=0$.
Moreover $A^*,A$ is isomorphic to (\ref{eq:affLP}) if and only if
$\xi= - \alpha$, $\zeta=\theta^*_d + \alpha \theta_0$, 
$\xi^*= - \alpha^{-1}$, $\zeta^* = \theta_0 + \alpha^{-1} \theta^*_d$.
\item[(iv)]
Assume $\phi_1 = \phi_d$ and $\varphi_1 \neq \varphi_d$.
Then $A,A^*$ is isomorphic to (\ref{eq:affLP}) if and only if
$\xi=1$, $\zeta= 0$, $\xi^*=1$, $\zeta^*=0$. 
Moreover $A^*,A$ is isomorphic to (\ref{eq:affLP}) if and only if
$\xi=\alpha$, $\zeta=\theta^*_0 - \alpha \theta_0$,  
$\xi^*=\alpha^{-1}$, $\zeta^*=\theta_0 - \alpha^{-1} \theta^*_0$.
\item[(v)]
Assume $\varphi_1=-\phi_1$, $\varphi_d=-\phi_d$ and $\varphi_1 \neq \varphi_d$.
Then $A,A^*$ is isomorphic to (\ref{eq:affLP}) if and only if
the sequence $\xi, \zeta, \xi^*, \zeta^*$  is listed in the following table.
\[
 \begin{array}{cccc}
  \;\;\; \xi \;\;\; & \zeta & \;\;\; \xi^* \;\;\; & \zeta^* \\
 \hline
   1 & 0 & 1 & 0 \\
  -1 & \theta_0+\theta_d & 1 & 0
 \end{array}
\]
Moreover $A^*,A$ is not isomorphic to (\ref{eq:affLP}) for any
$\xi, \zeta, \xi^*, \zeta^*$.
\item[(vi)]
Assume $\varphi_1=-\phi_d$, $\varphi_d=-\phi_1$ and $\varphi_1 \neq \varphi_d$.
Then $A,A^*$ is isomorphic to (\ref{eq:affLP}) if and only if
the sequence $\xi, \zeta, \xi^*, \zeta^*$  is listed in the following table.
\[
 \begin{array}{cccc}
  \;\;\; \xi \;\;\; & \zeta & \;\;\; \xi^* \;\;\; & \zeta^* \\
 \hline
   1 & 0 & 1 & 0 \\
  1 & 0 & -1 & \theta^*_0+\theta^*_d
 \end{array}
\]
Moreover $A^*,A$ is not isomorphic to (\ref{eq:affLP}) for any
$\xi, \zeta, \xi^*, \zeta^*$.
\item[(vii)]
Assume  none of (i)--(vi) hold above.
Then $A,A^*$ is isomorphic to (\ref{eq:affLP}) if and only if
$\xi=1$, $\zeta=0$, $\xi^*=1$, $\zeta^*=0$.
Moreover $A^*,A$ is not isomorphic to (\ref{eq:affLP}) for any
$\xi, \zeta, \xi^*, \zeta^*$.
\end{itemize}
\end{theorem}

\begin{proof}
Routine consequence of Propositions \ref{prop:LPaff} and \ref{prop:LPaffs}.
\end{proof}

\section{The parameter arrrays in closed form}

\indent
In \cite{NT:balanced} and \cite{T:array}
the parameter array of a Leonard system is given in closed form. 
For the rest of this paper we consider how the results of previous
sections look in terms of this form.

\begin{notation}          \label{notation3}
Let
$\Phi=(A;\{E_i\}_{i=0}^d; A^*; \{E^*_i\}_{i=0}^d)$
denote a Leonard system over $\mathbb{K}$ and let
$(\{\theta_i\}_{i=0}^d; \{\theta^*_i\}_{i=0}^d;
  \{\varphi_i\}_{i=1}^d$; $\{\phi_i\}_{i=1}^d)$
denote the corresponding parameter array.
We assume $d \geq 3$.
\end{notation}

\begin{notation}          \label{notation:cases}
Referring to Notation \ref{notation3},
let $\overline{\mathbb{K}}$ denote the algebraic closure of $\mathbb{K}$
and let $q$ denote a nonzero scalar in  $\overline{\mathbb{K}}$ such that
$q+q^{-1}+1$ is equal to the common value of (\ref{eq:indep}).
We consider the following types:
\begin{center}
\begin{tabular}{c|c}
type & description \\
\hline
I & $q \neq 1$, $q \neq -1$ \\
II & $q=1$, $\text{\rm Char}(\mathbb{K}) \neq 2$ \\
$\;\;$III$^+$ & $q=-1$, $\text{\rm Char}(\mathbb{K}) \neq 2$, $d$  even \\
$\;\;$III$^-$ &  $q=-1$, $\text{\rm Char}(\mathbb{K}) \neq 2$, $d$  odd \\
IV & $q=1$, $\text{\rm Char}(\mathbb{K}) = 2$
\end{tabular}
\end{center}
\end{notation}

\section{Type I: $q \neq 1$ and $q \neq -1$}

\begin{lemma}   \cite[Theorem 6.1]{NT:balanced}   \label{lem:I}     \samepage
Referring to Notation \ref{notation3}, assume $\Phi$ is
Type I. Then
there exists unique scalars $\eta$, $\mu$, $h$, $\eta^*$, $\mu^*$, 
$h^*$, $\tau$ in 
$\overline{\mathbb{K}}$ such that
\begin{eqnarray}
\theta_i &=& \eta + \mu q^i + h q^{d-i},  
   \label{eq:Ith} \\
\theta^*_i &=& \eta^* + \mu^* q^i + h^* q^{d-i}
   \label{eq:Iths}
\end{eqnarray}
for $0 \leq i \leq d$ and
\begin{eqnarray}
\varphi_i &=& (q^i-1)(q^{d-i+1}-1)(\tau -\mu \mu^* q^{i-1} - hh^* q^{d-i}),
   \label{eq:Ivphi}  \\
\phi_i &=& (q^i-1)(q^{d-i+1}-1)(\tau -h \mu^* q^{i-1} - \mu h^* q^{d-i})
   \label{eq:Iphi}
\end{eqnarray}
 for $1 \leq i \leq d$.
\end{lemma}

\begin{remark}         \label{rem:I}       \samepage
Referring to Lemma \ref{lem:I}, for $1 \leq i \leq d$ we have $q^i \neq 1$;
otherwise $\varphi_i=0$ by (\ref{eq:Ivphi}).
For $0 \leq i \leq d-1$ we have $\mu \neq h q^i$;
otherwise $\theta_{d-i}=\theta_0$. 
Similarly $\mu^* \neq h^* q^i$.
\end{remark}

\begin{lemma}              \label{lem:Itheta}   \samepage
Referring to Notation \ref{notation3}, assume $\Phi$ is
Type I. Then (i)--(iv) hold below.
\begin{itemize}
\item[(i)] $\theta_i + \theta_{d-i}$ is independent of $i$ for $0 \leq i \leq d$ 
if and only if $\mu=-h$.
\item[(ii)]
$\theta^*_i + \theta^*_{d-i}$ is independent of $i$ 
for $0 \leq i \leq d$ if and only if $\mu^*=-h^*$.
\item[(iii)]
$(\theta^*_i-\theta^*_0)(\theta_i-\theta_0)^{-1}$ is independent of $i$
for $1 \leq i \leq d$ if and only if $\mu h^*=\mu^*h$.
\item[(iv)]
$(\theta^*_{d-i}-\theta^*_d)(\theta_i-\theta_0)^{-1}$ is independent of $i$
for $1 \leq i \leq d$ if and only if $\mu\mu^*=hh^*$.
\end{itemize}
\end{lemma}

\begin{proof}
(i): Using (\ref{eq:Ith}),
\[
  \theta_i+\theta_{d-i}-\theta_0-\theta_d
  = (q^i-1)(1-q^{d-i})(\mu+h)
\]
for $0 \leq i \leq d$. 
The result follows from this and Remark \ref{rem:I}.

(ii): Similar to the proof of (i).

(iii):  Using (\ref{eq:Ith}) and (\ref{eq:Iths}),
\[
  \frac{\theta^*_i-\theta^*_0}
       {\theta_i-\theta_0}
 - \frac{\theta^*_d-\theta^*_0}
       {\theta_d-\theta_0}
 =\frac{(\mu h^*-\mu^*h)(1-q^{d-i})}
       {(\mu-h)(\mu-h q^{d-i})}
\]
for $1 \leq i \leq d$. 
The result follows from this and Remark \ref{rem:I}.

(iv): Using (\ref{eq:Ith}) and (\ref{eq:Iths}),
\[
  \frac{\theta^*_{d-i}-\theta^*_d}
       {\theta_i-\theta_0}
 - \frac{\theta^*_0-\theta^*_d}
       {\theta_d-\theta_0}
 =\frac{(\mu\mu^*-hh^*)(1-q^{d-i})}
       {(\mu-h)(\mu-h q^{d-i})}
\]
for $1 \leq i \leq d$. 
The result follows from this and Remark \ref{rem:I}.
\end{proof}

\begin{lemma}              \label{lem:Iphi}     \samepage
Referring to Notation \ref{notation3}, assume $\Phi$ is
Type I. Then (i)--(iv) hold below.
\begin{itemize}
\item[(i)]
$\varphi_i = - \phi_i$ for $1 \leq i \leq d$
if and only if $\tau=0$ and $\mu=-h$.
\item[(ii)]
$\varphi_i = - \phi_{d-i+1}$ for $1 \leq i \leq d$
if and only if $\tau=0$ and $\mu^*=-h^*$.
\item[(iii)]
$\phi_i = \phi_{d-i+1}$ for $1 \leq i \leq d$
if and only if $\mu h^*=\mu^*h$.
\item[(iv)]
$\varphi_i = \varphi_{d-i+1}$ for $1 \leq i \leq d$
if and only if $\mu\mu^*=hh^*$.
\end{itemize}
\end{lemma}

\begin{proof}
(i): 
Using the data in Lemma \ref{lem:I} we find
\begin{equation}          \label{eq:Iaux1}
 \varphi_i+\phi_i
 = (q^i-1)(q^{d-i+1}-1)(2\tau - (\mu+h)(\mu^* q^{i-1} + h^* q^{d-i}))
\end{equation}
for $1 \leq i \leq d$ and
\begin{equation}          \label{eq:Iaux2}
  \varphi_1+\phi_1-\varphi_d-\phi_d
 = (q-1)(q^{d-1}-1)(q^d-1)(\mu+h)(\mu^*-h^*).
\end{equation}
First assume $\varphi_i=-\phi_i$ for $1 \leq i \leq d$. 
Then in (\ref{eq:Iaux2}) the expression on the left is zero so the expression
on the right is zero. In this expression each factor except
$\mu+h$ is nonzero by Remark \ref{rem:I}, so $\mu=-h$.
In (\ref{eq:Iaux1}) the expression on the left is zero so the
expression on the right is zero. Evaluating this
expression using $\mu=-h$ and Remark \ref{rem:I} we find
$2 \tau = 0$. 
Note that $\text{Char}(\mathbb{K}) \neq 2$; 
otherwise $\mu=h$ and Remark \ref{rem:I} is contradicted. 
Therefore $\tau=0$.
We have now shown $\tau=0$ and $\mu=-h$. 
Conversely assume $\tau=0$ and $\mu=-h$. 
Then by (\ref{eq:Iaux1}) we have
$\varphi_i = - \phi_i$ for $1 \leq i \leq d$.

(ii):
Similar to the proof of (i). We note that
\[
\varphi_i+\phi_{d-i+1}
 = (q^i-1)(q^{d-i+1}-1)(2\tau - (\mu^*+h^*)(\mu q^{i-1} + h q^{d-i}))
\]
for $1 \leq i \leq d$ and
\[
\varphi_1+\phi_d - \varphi_d-\phi_1
=(q-1)(q^{d-1}-1)(q^d-1)(\mu-h)(\mu^*+h^*).
\]

(iii): 
Using the data in Lemma \ref{lem:I} we find
\[
 \phi_i-\phi_{d-i+1}
  = (q^i-1)(q^{d-i+1}-1)(q^{i-1}-q^{d-i})(\mu h^*-\mu^*h)
\]
for $1 \leq i \leq d$. 
The result follows from this and Remark \ref{rem:I}.

(iv):
Using the data in Lemma \ref{lem:I} we find
\[
 \varphi_i-\varphi_{d-i+1}
 = (q^i-1)(q^{d-i+1}-1)(q^{d-i}-q^{i-1})(\mu\mu^* - hh^*)
\]
for $1 \leq i \leq d$. 
The result follows from this and Remark \ref{rem:I}.
\end{proof}

\begin{proposition}        \label{prop:I}          \samepage
Referring to Notation \ref{notation3}, assume $\Phi$ is
Type I. Then (i)--(vii) hold below.
\begin{itemize}
\item[(i)]
$\Phi^\downarrow$ is affine isomorphic to $\Phi$ 
if and only if $\mu^*=-h^*$, $\tau=0$.
\item[(ii)]
$\Phi^\Downarrow$ is affine isomorphic to $\Phi$ 
if and only if $\mu=-h$, $\tau=0$.
\item[(iii)]
$\Phi^{\downarrow\Downarrow}$ is affine isomorphic to $\Phi$ 
if and only if  $\mu=-h$, $\mu^*=-h^*$.
\item[(iv)]
$\Phi^*$ is affine isomorphic to $\Phi$ 
if and only if $\mu h^*=\mu^*h$.
\item[(v)]
$\Phi^{\downarrow*}$ is affine isomorphic to $\Phi$ 
if and only if $\mu=-h$, $\mu^*=-h^*$, $\tau=0$.
\item[(vi)]
$\Phi^{\Downarrow*}$ is affine isomorphic to $\Phi$ 
if and only if $\mu=-h$, $\mu^*=-h^*$, $\tau=0$.
\item[(vii)]
$\Phi^{\downarrow\Downarrow*}$ is affine isomorphic to $\Phi$ 
if and only if $\mu\mu^*=hh^*$.
\end{itemize}
\end{proposition}

\begin{proof}
Follows from Propositions 
\ref{prop:down}, \ref{prop:Down}, \ref{prop:star}, \ref{prop:downDownstar},
\ref{prop:downDown}, \ref{prop:downstar}, \ref{prop:Downstar},
and Lemmas \ref{lem:Itheta}, \ref{lem:Iphi}.
\end{proof}

\begin{theorem}          \label{thm:I}   \samepage
Referring to Notation \ref{notation3}, assume $\Phi$ is
Type I.
Then (i)--(vii) hold below.
\begin{itemize}
\item[(i)]
Case (i) of Theorem \ref{thm:main} occurs if and only if
$\mu=-h$, $\mu^*=-h^*$ and $\tau=0$.
\item[(ii)]
Case (ii) of Theorem \ref{thm:main} occurs if and only if
$\mu=-h$, $\mu^*=-h^*$ and $\tau \neq 0$.
\item[(iii)]
Case (iii) of Theorem \ref{thm:main} occurs if and only if
$\mu\mu^*=hh^*$ and $\mu h^* \neq \mu^*h$.
\item[(iv)]
Case (iv) of Theorem \ref{thm:main} occurs if and only if
$\mu\mu^* \neq hh^*$ and $\mu h^*=\mu^*h$.
\item[(v)]
Case (v) of Theorem \ref{thm:main} occurs if and only if
$\mu=-h$, $\mu^*\neq -h^*$ and $\tau=0$.
\item[(vi)]
Case (vi) of Theorem \ref{thm:main} occurs if and only if
$\mu \neq -h$, $\mu^*=-h^*$ and $\tau=0$.
\item[(vii)]
Case (vii) of Theorem \ref{thm:main} occurs if and only if
$\mu\mu^* \neq hh^*$, $\mu h^* \neq \mu^*h$, and
at least two of $\mu\neq -h$, $\mu^* \neq -h^*$, $\tau\neq 0$.
\end{itemize}
\end{theorem}

\begin{proof}
Combine Theorem \ref{thm:main} and Proposition \ref{prop:I}.
\end{proof}

\section{Type II: $q = 1$ and $\text{Char}(\mathbb{K}) \neq 2$}

\begin{lemma}      \cite[Theorem 7.1]{NT:balanced}   \label{lem:II}  \samepage
Referring to Notation \ref{notation3}, assume $\Phi$ is
Type II. Then there exists unique scalars
$\eta$, $\mu$, $h$, $\eta^*$, $\mu^*$, $h^*$, $\tau$ in 
$\overline{\mathbb{K}}$ such that
\begin{eqnarray}
\theta_i &=& \eta + \mu(i-d/2)+ h i(d-i),
   \label{eq:IIth}  \\
\theta^*_i &=& \eta^* + \mu^*(i-d/2)+h^* i(d-i)
   \label{eq:IIths} 
\end{eqnarray}
for $0 \leq i \leq d$ and 
\begin{eqnarray}
\varphi_i &=& 
 i(d-i+1)(\tau-\mu \mu^*/2+(h \mu^*+ \mu h^*)(i-(d+1)/2)+hh^*(i-1)(d-i)),
   \label{eq:IIvphi}  \\
\phi_i &=& 
    i(d-i+1)(\tau+\mu \mu^*/2+(h\mu^*-\mu h^*)(i-(d+1)/2)+hh^*(i-1)(d-i))
   \label{eq:IIphi}
\end{eqnarray}
for $1 \leq i \leq d$.
\end{lemma}

\begin{remark}          \label{rem:II}     \samepage
Referring to Lemma \ref{lem:II}, 
for $0 \leq i \leq d-1$ we have $\mu \neq -ih$;
otherwise $\theta_{d-i}=\theta_0$.
Similarly $\mu^* \neq -i h^*$.
For any prime $i$ such that $i \leq d$ we have 
$\text{\rm Char}(\mathbb{K}) \neq i$; otherwise $\varphi_i=0$ by
(\ref{eq:IIvphi}).
\end{remark}

\begin{lemma}              \label{lem:IItheta}   \samepage
Referring to Notation \ref{notation3}, assume $\Phi$ is
Type II. Then (i)--(iv) hold below.
\begin{itemize}
\item[(i)] $\theta_i + \theta_{d-i}$ is independent of $i$ for $0 \leq i \leq d$ 
if and only if $h=0$.
\item[(ii)]
$\theta^*_i + \theta^*_{d-i}$ is independent of $i$ 
for $0 \leq i \leq d$ if and only if $h^*=0$.
\item[(iii)]
$(\theta^*_i-\theta^*_0)(\theta_i-\theta_0)^{-1}$ is independent of $i$
for $1 \leq i \leq d$ if and only if $\mu h^*=\mu^* h$.
\item[(iv)]
$(\theta^*_{d-i}-\theta^*_d)(\theta_i-\theta_0)^{-1}$ is independent of $i$
for $1 \leq i \leq d$ if and only if $\mu h^*=-\mu^* h$.
\end{itemize}
\end{lemma}

\begin{proof}
(i): Using (\ref{eq:IIth}),
\[
  \theta_i+\theta_{d-i}-\theta_0-\theta_d= 2i(d-i)h
\]
for $0 \leq i \leq d$. 
The result follows from this and Remark \ref{rem:II}.

(ii): Similar to the proof of (i).

(iii):  Using (\ref{eq:IIth}) and (\ref{eq:IIths}),
\[
  \frac{\theta^*_i-\theta^*_0}
       {\theta_i-\theta_0}
 - \frac{\theta^*_d-\theta^*_0}
       {\theta_d-\theta_0}
 =\frac{(d-i)(\mu h^*-\mu^*h)}
       {\mu(\mu+(d-i)h)}
\]
for $1 \leq i \leq d$. 
The result follows from this and Remark \ref{rem:II}.

(iv):  Using (\ref{eq:IIth}) and (\ref{eq:IIths}),
\[
  \frac{\theta^*_{d-i}-\theta^*_d}
       {\theta_i-\theta_0}
 - \frac{\theta^*_0-\theta^*_d}
       {\theta_d-\theta_0}
 =\frac{(d-i)(\mu h^* + \mu^*h)}
       {\mu(\mu+(d-i)h)}
\]
for $1 \leq i \leq d$. 
The result follows from this and Remark \ref{rem:II}.
\end{proof}

\begin{lemma}          \label{lem:IIphi}     \samepage
Referring to Notation \ref{notation3}, assume $\Phi$ is
Type II.
Then (i)--(iv) hold below.
\begin{itemize}
\item[(i)]
$\varphi_i = - \phi_i$ for $1 \leq i \leq d$
if and only if $\tau=0$ and $h=0$.
\item[(ii)]
$\varphi_i = - \phi_{d-i+1}$ for $1 \leq i \leq d$
if and only if $\tau=0$ and $h^*=0$.
\item[(iii)]
$\phi_i = \phi_{d-i+1}$ for $1 \leq i \leq d$
if and only if $\mu h^*=\mu^* h$.
\item[(iv)]
$\varphi_i = \varphi_{d-i+1}$ for $1 \leq i \leq d$
if and only if $\mu h^*=- \mu^* h$.
\end{itemize}
\end{lemma}

\begin{proof}
(i): 
Using the data in Lemma \ref{lem:II} we find
\begin{equation}          \label{eq:IIaux1}
 \varphi_i+\phi_i=2i(d-i+1)(\tau+\mu^*h(i-(d+1)/2)+hh^*(i-1)(d-i))
\end{equation}
for $1 \leq i \leq d$ and
\begin{equation}          \label{eq:IIaux2}
\varphi_1+\phi_1-\varphi_d-\phi_d=2d(1-d)\mu^*h.
\end{equation}
First assume $\varphi_i=-\phi_i$ for $1 \leq i \leq d$. 
Then in (\ref{eq:IIaux2}) the expression on the left is zero so the expression
on the right is zero. In this expression each factor except
$h$ is nonzero by Remark \ref{rem:II}, so $h=0$.
In (\ref{eq:IIaux1}) the expression on the left is zero so the
expression on the right is zero. Evaluating this
expression using $h=0$ and Remark \ref{rem:II} we find
$\tau=0$.
We have now shown $\tau=0$ and $h=0$. 
Conversely assume $\tau=0$ and $h=0$. 
Then by (\ref{eq:IIaux1}) we have
$\varphi_i = - \phi_i$ for $1 \leq i \leq d$.

(ii):
Similar to the proof of (i). We note that
\[
 \varphi_i+\phi_{d-i+1}=2i(d-i+1)(\tau+\mu h^*(i-(d+1)/2)+hh^*(i-1)(d-i))
\]
for $1 \leq i \leq d$ and
\[
 \varphi_1+\phi_d-\varphi_d-\phi_1=2d(1-d)\mu h^*.
\]

(iii):
Using the data in Lemma \ref{lem:II} we find
\[
 \phi_i-\phi_{d-i+1}=i(d-i+1)(d-2i+1)(\mu h^* - \mu^* h)
\]
for $1 \leq i \leq d$.
The result follows from this and Remark \ref{rem:II}.

(iv):
Using the data in Lemma \ref{lem:II} we find
\[
  \varphi_i-\varphi_{d-i+1}=-i(d-i+1)(d-2i+1)(\mu h^*+\mu^* h)
\]
for $1 \leq i \leq d$. 
The result follows from this and Remark \ref{rem:II}.
\end{proof}

\begin{proposition}         \label{prop:II}    \samepage
Referring to Notation \ref{notation3}, assume $\Phi$ is
Type II.
Then (i)--(vii) hold below.
\begin{itemize}
\item[(i)]
$\Phi^\downarrow$ is affine isomorphic to $\Phi$ 
if and only if $h^*=0$, $\tau=0$.
\item[(ii)]
$\Phi^\Downarrow$ is affine isomorphic to $\Phi$ 
if and only if $h=0$, $\tau=0$.
\item[(iii)]
$\Phi^{\downarrow\Downarrow}$ is affine isomorphic to $\Phi$ 
if and only if  $h=0$, $h^*=0$.
\item[(iv)]
$\Phi^*$ is affine isomorphic to $\Phi$ 
if and only if $\mu h^*= \mu^* h$.
\item[(v)]
$\Phi^{\downarrow*}$ is affine isomorphic to $\Phi$ 
if and only if $h=0$, $h^*=0$, $\tau=0$.
\item[(vi)]
$\Phi^{\Downarrow*}$ is affine isomorphic to $\Phi$ 
if and only if $h=0$, $h^*=0$, $\tau=0$.
\item[(vii)]
$\Phi^{\downarrow\Downarrow*}$ is affine isomorphic to $\Phi$ 
if and only if $\mu h^*=-\mu^* h$.
\end{itemize}
\end{proposition}

\begin{proof}
Follows from Propositions 
\ref{prop:down}, \ref{prop:Down}, \ref{prop:star}, \ref{prop:downDownstar},
\ref{prop:downDown}, \ref{prop:downstar}, \ref{prop:Downstar},
and Lemmas \ref{lem:IItheta}, \ref{lem:IIphi}.
\end{proof}

\begin{theorem}          \label{thm:II}   \samepage
Referring to Notation \ref{notation3}, assume $\Phi$ is
Type II.
Then (i)--(vii) hold below.
\begin{itemize}
\item[(i)]
Case (i) of Theorem \ref{thm:main} occurs if and only if
$h=0$, $h^*=0$ and $\tau=0$.
\item[(ii)]
Case (ii) of Theorem \ref{thm:main} occurs if and only if
$h=0$, $h^*=0$ and $\tau \neq 0$.
\item[(iii)]
Case (iii) of Theorem \ref{thm:main} occurs if and only if
$\mu h^* \neq \mu^* h$ and $\mu h^* = -\mu^* h$.
\item[(iv)]
Case (iv) of Theorem \ref{thm:main} occurs if and only if
$\mu h^* = \mu^* h$ and $\mu h^* \neq - \mu^* h$.
\item[(v)]
Case (v) of Theorem \ref{thm:main} occurs if and only if
$h=0$, $h^* \neq 0$ and $\tau=0$.
\item[(vi)]
Case (vi) of Theorem \ref{thm:main} occurs if and only if
$h \neq 0$, $h^* = 0$ and $\tau=0$.
\item[(vii)]
Case (vii) of Theorem \ref{thm:main} occurs if and only if
$\mu h^* \neq \mu^* h$, $\mu h^* \neq - \mu^* h$, and
at least two of $h$, $h^*$, $\tau$ are nonzero.
\end{itemize}
\end{theorem}

\begin{proof}
Combine Theorem \ref{thm:main} and Proposition \ref{prop:II}.
\end{proof}

\section{Type III$^+$: $q = -1$, $\text{Char}(\mathbb{K}) \neq 2$ and $d$ is even}

\begin{lemma}     \cite[Theorem 8.1]{NT:balanced}   \label{lem:IIIeven} \samepage
Referring to Notation \ref{notation3}, assume $\Phi$ is
Type III$^+$. Then there exists unique scalars 
$\eta$, $h$, $s$, $\eta^*$, $h^*$, $s^*$,
$\tau$ in $\overline{\mathbb{K}}$ such that
\begin{eqnarray}
\theta_i &=&
  \begin{cases}
     \eta+s+h(i-d/2)  & \text{\rm if $i$ is even}, \\
     \eta-s -h(i-d/2) & \text{\rm if $i$ is odd},
  \end{cases}   
         \label{eq:IIIeventh}   \\
\theta^*_i &=&
   \begin{cases}
     \eta^* +s^* +h^*(i-d/2)  &   \text{\rm if $i$ is even}, \\
     \eta^* - s^* -h^*(i-d/2) &  \text{\rm if $i$ is odd}
   \end{cases}
      \label{eq:IIIevenths}
\end{eqnarray}
for $0 \leq i \leq d$ and
\begin{eqnarray}    
\varphi_i &=&
   \begin{cases}
      i(\tau-s h^*-s^*h-hh^*(i-(d+1)/2)) &   \text{\rm if $i$ is even}, \\
      (d-i+1)(\tau+sh^*+s^*h+hh^*(i-(d+1)/2))  & \text{\rm if $i$ is odd},
   \end{cases}
     \label{eq:IIIevenvphi}  \\
\phi_i &=&
   \begin{cases}
      i(\tau-sh^*+s^*h+hh^*(i-(d+1)/2)) &   \text{\rm if $i$ is even}, \\
      (d-i+1)(\tau+sh^*-s^*h - hh^*(i-(d+1)/2))  & \text{\rm if $i$ is odd}
   \end{cases}
    \label{eq:IIIevenphi}
\end{eqnarray}
for $1 \leq i \leq d$.
\end{lemma}

\begin{remark}        \label{rem:IIIeven}     \samepage
Referring to Lemma \ref{lem:IIIeven},
we have $h \neq 0$; otherwise $\theta_0=\theta_2$ by (\ref{eq:IIIeventh}).
Similary we have $h^* \neq 0$.
For $i$ odd with $0 \leq i \leq d-1$ we have $s \neq i h/2$;
otherwise $\theta_{d-i}=\theta_0$.
For any prime $i$ such that $i \leq d/2$ we have 
$\text{\rm Char}(\mathbb{K})\neq i$; otherwise $\varphi_{2i}=0$
by ($\ref{eq:IIIevenvphi}$). By this and since $\text{\rm Char}(\mathbb{K})\neq 2$
we find $\text{\rm Char}(\mathbb{K})$ is either $0$ or an odd prime
greater than $d/2$.
Observe neither of $d$, $d-2$ vanish in $\mathbb{K}$ since otherwise
$\text{\rm Char}(\mathbb{K})$ must divide $d/2$ or $(d-2)/2$.
\end{remark}

\begin{lemma}              \label{lem:IIIeventheta}   \samepage
Referring to Notation \ref{notation3}, assume $\Phi$ is
Type III$^+$. Then (i)--(iv) hold below.
\begin{itemize}
\item[(i)] $\theta_i + \theta_{d-i}$ is independent of $i$ for $0 \leq i \leq d$ 
if and only if $s=0$.
\item[(ii)]
$\theta^*_i + \theta^*_{d-i}$ is independent of $i$ 
for $0 \leq i \leq d$ if and only if $s^*=0$.
\item[(iii)]
$(\theta^*_i-\theta^*_0)(\theta_i-\theta_0)^{-1}$ is independent of $i$
for $1 \leq i \leq d$ if and only if $h s^*=h^* s$.
\item[(iv)]
$(\theta^*_{d-i}-\theta^*_d)(\theta_i-\theta_0)^{-1}$ is independent of $i$
for $1 \leq i \leq d$ if and only if $h s^*=-h^* s$.
\end{itemize}
\end{lemma}

\begin{proof}
(i):
 Using (\ref{eq:IIIeventh}),
\[
  \theta_i + \theta_{d-i} =
   \begin{cases}
     2(\eta+s) & \text{ if $i$ is even}, \\
     2(\eta-s) & \text{ if $i$ is odd}
   \end{cases}
\]
for $0 \leq i \leq d$. 
The result follows from this.

(ii): Similar to the proof of (i).

(iii):
Using (\ref{eq:IIIeventh}) and (\ref{eq:IIIevenths}),
\[
  \frac{\theta^*_i - \theta^*_0}
       {\theta_i - \theta_0} 
 -\frac{\theta^*_d - \theta^*_0}
       {\theta_d-\theta_0}
 = \begin{cases}
     0 & \text{ if $i$ is even},  \\
     \displaystyle \frac{h s^* - h^* s}
          {h (s -h(d-i)/2)}       & \text{ if $i$ is odd}
   \end{cases}
\]
for $1 \leq i \leq d$. 
The result follows from this.

(iv):
Using (\ref{eq:IIIeventh}) and (\ref{eq:IIIevenths}),
\[
  \frac{\theta^*_{d-i}-\theta^*_d}
       {\theta_i-\theta_0}
 - \frac{\theta^*_0 - \theta^*_d}
        {\theta_d -\theta_0}
 = \begin{cases}
      0 & \text{ if $i$ is even},  \\
   \displaystyle \frac{h s^*+h^* s}
         {h(s- h(d-i)/2)}  &  \text{ if $i$ is odd}
   \end{cases}
\]
for $1 \leq i \leq d$. 
The result follows from this.
\end{proof}

\begin{lemma}              \label{lem:IIIevenphi}     \samepage
Referring to Notation \ref{notation3}, assume $\Phi$ is
Type III$^+$.
Then (i)--(iv) hold below.
\begin{itemize}
\item[(i)]
$\varphi_i = - \phi_i$ for $1 \leq i \leq d$
if and only if $\tau=0$ and $s=0$.
\item[(ii)]
$\varphi_i = - \phi_{d-i+1}$ for $1 \leq i \leq d$
if and only if $\tau=0$ and $s^*=0$.
\item[(iii)]
$\phi_i = \phi_{d-i+1}$ for $1 \leq i \leq d$
if and only if $h s^*=h^* s$.
\item[(iv)]
$\varphi_i = \varphi_{d-i+1}$ for $1 \leq i \leq d$
if and only if $h s^*=- h^* s$.
\end{itemize}
\end{lemma}

\begin{proof}
(i):
Using the data in Lemma \ref{lem:IIIeven} we find
\[
 \varphi_i+\phi_i =
 \begin{cases}
   2i(\tau-h^*s) & \text{ if $i$ is even}, \\
  2(d-i+1)(\tau+h^*s) & \text{ if $i$ is odd}
 \end{cases}
\]
for $1 \leq i \leq d$. 
The result follows from this and Remark \ref{rem:IIIeven}.

(ii):
Using the data in Lemma \ref{lem:IIIeven} we find
\[
 \varphi_i+\phi_{d-i+1} =
 \begin{cases}
   2i(\tau-h s^*) & \text{ if $i$ is even}, \\
  2(d-i+1)(\tau+h s^*) & \text{ if $i$ is odd}
 \end{cases}
\]
for $1 \leq i \leq d$. 
The result follows from this and Remark \ref{rem:IIIeven}.

(iii):
Using the data in Lemma \ref{lem:IIIeven} we find
\[
 \phi_i - \phi_{d-i+1} =
  \begin{cases}
    2i(h s^* - h^* s) & \text{ if $i$ is even}, \\
    -2(d-i+1)(h s^* - h^* s) & \text{ if $i$ is odd}
  \end{cases}
\]
for $1 \leq i \leq d$. 
The result follows from this and Remark \ref{rem:IIIeven}.

(iv):
Using the data in Lemma \ref{lem:IIIeven} we find
\[
 \varphi_i - \varphi_{d-i+1} =
  \begin{cases}
    -2i(h s^* + h^* s) & \text{ if $i$ is even}, \\
    2(d-i+1)(h s^* + h^* s) & \text{ if $i$ is odd}
  \end{cases}
\]
for $1 \leq i \leq d$. 
The result follows from this and Remark \ref{rem:IIIeven}.
\end{proof}

\begin{proposition}         \label{prop:IIIeven}         \samepage
Referring to Notation \ref{notation3}, assume $\Phi$ is
Type III$^+$.
Then (i)--(vii) hold below.
\begin{itemize}
\item[(i)]
$\Phi^\downarrow$ is affine isomorphic to $\Phi$ 
if and only if $s^*=0$, $\tau=0$.
\item[(ii)]
$\Phi^\Downarrow$ is affine isomorphic to $\Phi$ 
if and only if $s=0$, $\tau=0$.
\item[(iii)]
$\Phi^{\downarrow\Downarrow}$ is affine isomorphic to $\Phi$ 
if and only if  $s=0$, $s^*=0$.
\item[(iv)]
$\Phi^*$ is affine isomorphic to $\Phi$ 
if and only if $h s^*= h^* s$.
\item[(v)]
$\Phi^{\downarrow*}$ is affine isomorphic to $\Phi$ 
if and only if $s=0$, $s^*=0$, $\tau=0$.
\item[(vi)]
$\Phi^{\Downarrow*}$ is affine isomorphic to $\Phi$ 
if and only if $s=0$, $s^*=0$, $\tau=0$.
\item[(vii)]
$\Phi^{\downarrow\Downarrow*}$ is affine isomorphic to $\Phi$ 
if and only if $h s^*=-h^* s$.
\end{itemize}
\end{proposition}

\begin{proof}
Follows from Propositions 
\ref{prop:down}, \ref{prop:Down}, \ref{prop:star}, \ref{prop:downDownstar},
\ref{prop:downDown}, \ref{prop:downstar}, \ref{prop:Downstar},
and Lemmas \ref{lem:IIIeventheta}, \ref{lem:IIIevenphi}.
\end{proof}

\begin{theorem}          \label{thm:IIIeven}   \samepage
Referring to Notation \ref{notation3}, assume $\Phi$ is
Type III$^+$.
Then (i)--(vii) hold below.
\begin{itemize}
\item[(i)]
Case (i) of Theorem \ref{thm:main} occurs if and only if
$s=0$, $s^*=0$ and $\tau=0$.
\item[(ii)]
Case (ii) of Theorem \ref{thm:main} occurs if and only if
$s=0$, $s^*=0$ and $\tau \neq 0$.
\item[(iii)]
Case (iii) of Theorem \ref{thm:main} occurs if and only if
$h s^* \neq h^* s$ and $h s^* = -h^* s$.
\item[(iv)]
Case (iv) of Theorem \ref{thm:main} occurs if and only if
$h s^* = h^* s$ and $h s^* \neq - h^* s$.
\item[(v)]
Case (v) of Theorem \ref{thm:main} occurs if and only if
$s=0$, $s^* \neq 0$ and $\tau=0$.
\item[(vi)]
Case (vi) of Theorem \ref{thm:main} occurs if and only if
$s \neq 0$, $s^* = 0$ and $\tau=0$.
\item[(vii)]
Case (vii) of Theorem \ref{thm:main} occurs if and only if
$h s^* \neq h^* s$, $h s^* \neq - h^* s$, and
at least two of $s$, $s^*$, $\tau$ are nonzero.
\end{itemize}
\end{theorem}

\begin{proof}
Combine Theorem \ref{thm:main} and Proposition \ref{prop:IIIeven}.
\end{proof}

\section{Type III$^-$: $q = -1$, $\text{Char}(\mathbb{K}) \neq 2$ and $d$ is odd}

\begin{lemma}   \cite[Theorem 9.1]{NT:balanced}   \label{lem:IIIodd}  \samepage
Referring to Notation \ref{notation3}, assume $\Phi$ is
Type III$^-$.
Then there exists unique scalars 
$\eta$, $h$, $s$, $\eta^*$, $h^*$, $s^*$, $\tau$
in $\overline{\mathbb{K}}$ such that 
\begin{eqnarray}
\theta_i &=&
  \begin{cases}
     \eta+s+h(i-d/2)  & \text{\rm if $i$ is even}, \\
     \eta-s-h(i-d/2) & \text{\rm if $i$ is odd},
  \end{cases}   
         \label{eq:IIIoddth}   \\
\theta^*_i &=&
   \begin{cases}
     \eta^* + s^* +h^*(i-d/2)  &   \text{\rm if $i$ is even}, \\
     \eta^*-s^* -h^*(i-d/2) &  \text{\rm if $i$ is odd}
   \end{cases}
      \label{eq:IIIoddths}
\end{eqnarray}
for $0 \leq i \leq d$ and 
\begin{eqnarray}    
\varphi_i &=&
   \begin{cases}
      hh^* i(d-i+1) &   \text{\rm if $i$ is even}, \\
      \tau-2ss^*+i(d-i+1)hh^*-2(h s^*+h^*s)(i-(d+1)/2)  & \text{\rm if $i$ is odd},
   \end{cases}
      \label{eq:IIIoddvphi}  \\
\phi_i &=&
   \begin{cases}
      hh^* i(d-i+1) &   \text{\rm if $i$ is even}, \\
      \tau+2ss^*+i(d-i+1)hh^*-2(h s^*-h^*s)(i-(d+1)/2)  & \text{\rm if $i$ is odd}
   \end{cases}
      \label{eq:IIIoddphi}
\end{eqnarray}
for $1 \leq i \leq d$.
\end{lemma}

\begin{remark}         \label{rem:IIIodd}     \samepage
Referring to Lemma \ref{lem:IIIodd},
we have $h \neq 0$; otherwise $\theta_0=\theta_2$ by (\ref{eq:IIIoddth}).
Similary we have $h^* \neq 0$.
We have $s \neq 0$; otherwise $\theta_0=\theta_d$ by (\ref{eq:IIIoddth}).
Similarly we have $s^* \neq 0$.
For $i$ even with $0 \leq i \leq d-1$ we have $s \neq ih/2$;
otherwise $\theta_{d-i}=\theta_0$.
For any prime $i$ such that $i \leq d/2$ we have 
$\text{\rm Char}(\mathbb{K})\neq i$; otherwise $\varphi_{2i}=0$
by ($\ref{eq:IIIoddvphi}$). By this and since $\text{\rm Char}(\mathbb{K})\neq 2$
we find $\text{\rm Char}(\mathbb{K})$ is either $0$ or an odd prime
greater than $d/2$.
Observe $d-1$ does not vanish in $\mathbb{K}$ since otherwise
$\text{\rm Char}(\mathbb{K})$ must divide $(d-1)/2$.
\end{remark}

\begin{lemma}              \label{lem:IIIoddtheta}   \samepage
Referring to Notation \ref{notation3}, assume $\Phi$ is
Type III$^-$. Then (i)--(iv) hold below.
\begin{itemize}
\item[(i)] 
$\theta_0+\theta_d \neq \theta_1+\theta_{d-1}$.
\item[(ii)]
$\theta^*_0 + \theta^*_d \neq \theta^*_1 +\theta^*_{d-1}$.
\item[(iii)]
$(\theta^*_i-\theta^*_0)(\theta_i-\theta_0)^{-1}$ is independent of $i$
for $1 \leq i \leq d$ if and only if $h s^*=h^* s$.
\item[(iv)]
$(\theta^*_{d-i}-\theta^*_d)(\theta_i-\theta_0)^{-1}$ is independent of $i$
for $1 \leq i \leq d$ if and only if $h s^*=-h^* s$.
\end{itemize}
\end{lemma}

\begin{proof}
(i):
 Using (\ref{eq:IIIoddth}),
\[
 \theta_0+\theta_d-\theta_1-\theta_{d-1}=-2(d-1)h.
\]
The result follows from this and Remark \ref{rem:IIIodd}.

(ii): Similar to the proof of (i).

(iii):
Using (\ref{eq:IIIoddth}) and (\ref{eq:IIIoddths}),
\[
  \frac{\theta^*_i - \theta^*_0}
       {\theta_i - \theta_0} 
 -\frac{\theta^*_d - \theta^*_0}
       {\theta_d-\theta_0}
 = \begin{cases}
     \displaystyle 
     \frac{h^* s - h s^*}
          {hs}                    & \text{ if $i$ is even},  \\
    \ \\
     \displaystyle 
     \frac{(d-i)(h s^* - h^* s)}
          {s (2s -h(d-i))}       & \text{ if $i$ is odd}
   \end{cases}
\]
for $1 \leq i \leq d$. 
The result follows from this.

(iv):
Using (\ref{eq:IIIoddth}) and (\ref{eq:IIIoddths}),
\[
  \frac{\theta^*_{d-i}-\theta^*_d}
       {\theta_i-\theta_0}
 - \frac{\theta^*_0 - \theta^*_d}
        {\theta_d -\theta_0}
 = \begin{cases}
   \displaystyle \frac{hs^*+h^* s}{hs}  & \text{ if $i$ is even},  \\
   \\
   \displaystyle 
       \frac{(d-i)(hs^*+h^*s)}{s(-2s+h(d-i))}  &  \text{ if $i$ is odd}
   \end{cases}
\]
for $1 \leq i \leq d$. 
The result follows from this.
\end{proof}

\begin{lemma}      \label{lem:IIIoddphi}     \samepage
Referring to Notation \ref{notation3}, assume $\Phi$ is
Type III$^-$. Then (i)--(iv) hold below.
\begin{itemize}
\item[(i)] 
$\varphi_2 \neq - \phi_2$.
Moreover if $\varphi_1 = - \phi_1$ then $\varphi_d \neq - \phi_d$.
\item[(ii)]
$\varphi_2 \neq - \phi_{d-1}$.
Moreover if  $\varphi_1 = - \phi_d$ then $\varphi_d \neq - \phi_1$.
\item[(iii)]
$\phi_i = \phi_{d-i+1}$ for $1 \leq i \leq d$
if and only if $h s^*=h^* s$.
\item[(iv)]
$\varphi_i = \varphi_{d-i+1}$ for $1 \leq i \leq d$
if and only if $h s^*=- h^* s$.
\end{itemize}
\end{lemma}

\begin{proof}
(i):
Using the data in Lemma \ref{lem:IIIodd} we find
\[
  \varphi_2 + \phi_2 = 4(d-1)hh^*,
\]
\[
  \varphi_1+ \phi_1-\varphi_d-\phi_d
 = 4(d-1)h s^*.
\]
The result follows from this and Remark \ref{rem:IIIodd}.

(ii):
Using the data in Lemma \ref{lem:IIIodd} we find
\[
  \varphi_2 + \phi_{d-1} = 4(d-1)hh^*,
\]
\[
  \varphi_1+\phi_d-\varphi_d-\phi_1
 = 4(d-1)h^* s.
\]
The result follows from this and Remark \ref{rem:IIIodd}.

(iii):
Using the data in Lemma \ref{lem:IIIodd} we find
\[
\phi_i - \phi_{d-i+1} =
 \begin{cases}
   0 & \text{ if $i$ is even},   \\
  2(d-2i+1)(h s^* - h^* s) & \text{ if $i$ is odd}
 \end{cases}
\]
for $1 \leq i \leq d$. 
The result follows from this and Remark \ref{rem:IIIodd}.

(iv):
Using the data in Lemma \ref{lem:IIIodd} we find
\[
\varphi_i - \varphi_{d-i+1} =
 \begin{cases}
   0 & \text{ if $i$ is even},   \\
  2(d-2i+1)(h s^* + h^* s) & \text{ if $i$ is odd}
 \end{cases}
\]
for $1 \leq i \leq d$. 
The result follows from this and Remark \ref{rem:IIIodd}.
\end{proof}

\begin{proposition}         \label{prop:IIIodd}          \samepage
Referring to Notation \ref{notation3}, assume $\Phi$ is
Type III$^-$.
Then (i)--(vii) hold below.
\begin{itemize}
\item[(i)]
$\Phi^\downarrow$ is not affine isomorphic to $\Phi$.
\item[(ii)]
$\Phi^\Downarrow$ is not affine isomorphic to $\Phi$. 
\item[(iii)]
$\Phi^{\downarrow\Downarrow}$ is not affine isomorphic to $\Phi$.
\item[(iv)]
$\Phi^*$ is affine isomorphic to $\Phi$ 
if and only if $h s^*= h^* s$.
\item[(v)]
$\Phi^{\downarrow*}$ is not affine isomorphic to $\Phi$.
\item[(vi)]
$\Phi^{\Downarrow*}$ is not affine isomorphic to $\Phi$.
\item[(vii)]
$\Phi^{\downarrow\Downarrow*}$ is affine isomorphic to $\Phi$ 
if and only if $h s^*=-h^* s$.
\end{itemize}
\end{proposition}

\begin{proof}
Follows from Propositions 
\ref{prop:down}, \ref{prop:Down}, \ref{prop:star}, \ref{prop:downDownstar},
\ref{prop:downDown}, \ref{prop:downstar}, \ref{prop:Downstar},
and Lemmas \ref{lem:IIIoddtheta}, \ref{lem:IIIoddphi}.
\end{proof}

\begin{theorem}          \label{thm:IIIodd}   \samepage
Referring to Notation \ref{notation3}, assume $\Phi$ is
Type III$^-$.
Then (i)--(iv) hold below.
\begin{itemize}
\item[(i)]
Case (iii) of Theorem \ref{thm:main} occurs if and only if
$h s^* = - h^* s$.
\item[(ii)]
Case (iv) of Theorem \ref{thm:main} occurs if and only if
$h s^* = h^* s$.
\item[(iii)]
Case (vii) of Theorem \ref{thm:main} occurs if and only if
both $h s^* \neq h^* s$, $h s^* \neq - h^* s$.
\item[(iv)]
Cases (i), (ii), (v), (vi) of  Theorem \ref{thm:main} do not occur.
\end{itemize}
\end{theorem}

\begin{proof}
Combine Theorem \ref{thm:main} and Proposition \ref{prop:IIIodd}.
\end{proof}

\section{Type IV: $q = 1$ and $\text{Char}(\mathbb{K}) = 2$}

\begin{lemma}  \cite[Theorem 10.1]{NT:balanced}   \label{lem:IV}     \samepage
Referring to Notation \ref{notation3}, assume $\Phi$ is
Type IV.
Then $d=3$. Moreover there exists unique scalars
$h$, $s$, $h^*$, $s^*$, $r$ in $\overline{\mathbb{K}}$ such that
\[
 \begin{array}{lll}
    \theta_1 = \theta_0 + h(s+1),
  & \theta_2 = \theta_0 + h,
  & \theta_3 = \theta_0 + h s,
  \\
    \theta^*_1 = \theta^*_0 + h^*(s^*+1), 
  & \theta^*_2 = \theta^*_0 + h^*,
  & \theta^*_3 = \theta^*_0 + h^*s^*,
  \\
    \varphi_1 = hh^* r,
  & \varphi_2 = hh^*, 
  & \varphi_3 = hh^*(r+s+s^*), 
  \\
    \phi_1 = hh^*(r+s(1+s^*)),
  & \phi_2 = hh^*,
  & \phi_3 = hh^*(r+s^*(1+s)).
  \end{array}
\]
\end{lemma}

\begin{remark}     \label{rem:IV}       \samepage
Referring to Lemma \ref{lem:IV},
each of $h$, $h^*$, $s$, $s^*$ is nonzero, and each of
$s$, $s^*$ is not equal to $1$.
\end{remark}

\begin{lemma}              \label{lem:IVtheta}   \samepage
Referring to Notation \ref{notation3}, assume $\Phi$ is
Type IV. Then  (i)--(iv) hold below.
\begin{itemize}
\item[(i)] $\theta_i + \theta_{d-i}=h s$ for $0 \leq i \leq d$.
\item[(ii)]
$\theta^*_i + \theta^*_{d-i}=h^* s^*$ for $0 \leq i \leq d$.
\item[(iii)]
$(\theta^*_i-\theta^*_0)(\theta_i-\theta_0)^{-1}$ is independent of $i$
for $1 \leq i \leq d$ if and only if $s=s^*$.
\item[(iv)]
$(\theta^*_{d-i}-\theta^*_d)(\theta_i-\theta_0)^{-1}$ is independent of $i$
for $1 \leq i \leq d$ if and only if $s=s^*$.
\end{itemize}
\end{lemma}

\begin{proof}
(i), (ii): Routine verification using the data in Lemma \ref{lem:IV}.

(iii):
Using the data in Lemma \ref{lem:IV} we find
\[
\frac{\theta^*_1-\theta^*_0}
     {\theta_1-\theta_0} 
= \frac{h^*(s^*+1)}
       {h(s+1)},
\qquad
\frac{\theta^*_2-\theta^*_0}
     {\theta_2-\theta_0}
= \frac{h^*}
       {h},
\qquad
\frac{\theta^*_3-\theta^*_0}
     {\theta_3-\theta_0}
= \frac{h^* s^*}
       {h s}.
\]
The result follows from this and Remark \ref{rem:IV}.

(iv):
Using the data in Lemma \ref{lem:IV} we find
\[
\frac{\theta^*_2-\theta^*_3}
     {\theta_1-\theta_0} 
= \frac{h^*(s^*+1)}
       {h(s+1)},
\qquad
\frac{\theta^*_1-\theta^*_3}
     {\theta_2-\theta_0}
= \frac{h^*}
       {h},
\qquad
\frac{\theta^*_0-\theta^*_3}
     {\theta_3-\theta_0}
= \frac{h^* s^*}
       {h s}.
\]
The result follows from this and Remark \ref{rem:IV}.
\end{proof}

\begin{lemma}         \label{lem:IVphi}     \samepage
Referring to Notation \ref{notation3}, assume $\Phi$ is
Type IV. Then (i)--(iv) hold below.
\begin{itemize}
\item[(i)]
$\varphi_1 \neq - \phi_1$.
\item[(ii)]
$\varphi_1 \neq - \phi_3$.
\item[(iii)]
$\phi_1 = \phi_3$ if and only if $s=s^*$.
\item[(iv)]
$\varphi_1 = \varphi_3$ if and only if $s=s^*$.
\end{itemize}
\end{lemma}

\begin{proof}
Using the data in Lemma \ref{lem:IV} we find
\[
 \varphi_1+\phi_1=hh^*s(s^*+1),
 \qquad
 \varphi_1+\phi_3=hh^*s^*(s+1),
\]
\[
\phi_1-\phi_3=hh^*(s+s^*),
\qquad
\varphi_1-\varphi_3=hh^*(s+s^*).
\]
Now (i)--(iv) follow from this and Remark \ref{rem:IV}.
\end{proof}

\begin{proposition}         \label{prop:IV}         \samepage
Referring to Notation \ref{notation3}, assume $\Phi$ is
Type IV.
Then (i)--(vii) hold below.
\begin{itemize}
\item[(i)]
$\Phi^\downarrow$ is not affine isomorphic to $\Phi$.
\item[(ii)]
$\Phi^\Downarrow$ is not affine isomorphic to $\Phi$. 
\item[(iii)]
$\Phi^{\downarrow\Downarrow}$ is affine isomorphic to $\Phi$ 
if and only if  $s=s^*$.
\item[(iv)]
$\Phi^*$ is affine isomorphic to $\Phi$ 
if and only if $s=s^*$.
\item[(v)]
$\Phi^{\downarrow*}$ is not affine isomorphic to $\Phi$.
\item[(vi)]
$\Phi^{\downarrow*}$ is not affine isomorphic to $\Phi$.
\item[(vii)]
$\Phi^{\downarrow\Downarrow*}$ is affine isomorphic to $\Phi$ 
if and only if $s=s^*$.
\end{itemize}
\end{proposition}

\begin{proof}
Follows from Propositions 
\ref{prop:down}, \ref{prop:Down}, \ref{prop:star}, \ref{prop:downDownstar},
\ref{prop:downDown}, \ref{prop:downstar}, \ref{prop:Downstar},
and Lemmas \ref{lem:IVtheta}, \ref{lem:IVphi}.
\end{proof}

\begin{theorem}          \label{thm:IV}   \samepage
Referring to Notation \ref{notation3}, assume $\Phi$ is
Type IV.
Then (i)--(iii) hold below.
\begin{itemize}
\item[(i)] 
Case (ii) of Theorem \ref{thm:main} occurs if and only if
$s=s^*$.
\item[(ii)] 
Case (vii) of Theorem \ref{thm:main} occurs if and only if
$s \neq s^*$.
\item[(iii)]
Cases (i), (iii)--(vi) of Theorem \ref{thm:main} do not occur.
\end{itemize}
\end{theorem}

\begin{proof}
Combine Theorem \ref{thm:main} and Proposition \ref{prop:IV}.
\end{proof}

\bigskip

\bibliographystyle{plain}

\bigskip\bigskip\noindent
Kazumasa Nomura\\
College of Liberal Arts and Sciences\\
Tokyo Medical and Dental University\\
Kohnodai, Ichikawa, 272-0827 Japan\\
email: knomura@pop11.odn.ne.jp

\bigskip\noindent
Paul Terwilliger\\
Department of Mathematics\\
University of Wisconsin\\
480 Lincoln Drive\\ 
Madison, Wisconsin, 53706 USA\\
email: terwilli@math.wisc.edu

\bigskip\noindent
{\bf Keywords.}
Leonard pair, tridiagonal pair, $q$-Racah polynomial, orthogonal polynomial.

\noindent
{\bf 2000 Mathematics Subject Classification}.
05E35, 05E30, 33C45, 33D45.

\end{document}